\RequirePackage{fix-cm}
\documentclass[journal,transmag,10pt]{IEEEtran}

\usepackage{listings}
\usepackage{color}
\usepackage{tabularx}
\definecolor{mygreen}{RGB}{28,172,0}
\definecolor{mylilas}{RGB}{170,55,241}
\usepackage{graphicx}
\usepackage{subfigure}
\usepackage{multirow}
\usepackage{url}
\usepackage{amssymb}
\usepackage{amsmath}
\usepackage{amsfonts}
\usepackage{wrapfig}
\usepackage{afterpage}
\usepackage{float}
\usepackage{algorithm}
\usepackage[noend]{algpseudocode}
\usepackage{epsfig}
\usepackage{epstopdf}

\usepackage{url}

\usepackage{listings}
\usepackage{color}
\definecolor{mygreen}{RGB}{28,172,0}
\definecolor{mylilas}{RGB}{170,55,241}

\graphicspath{{fg/}}
\newcommand{\real}{\mathbb{R}}
\newcommand{\beps}{\mbox{\boldmath$\epsilon$}}
\newcommand{\brho}{{\bf z}} 

\newcommand{\bphi}{\mbox{\boldmath$\phi$}}
\newcommand{\balp}{\mbox{\boldmath$\alpha$}}
\newcommand{\btau}{\mbox{\boldmath$\psi$}}
\def\vsig{\varsigma}
\newcommand{\bchi}{\mbox{\boldmath$\chi$}}
\newcommand{\bveps}{\mbox{\boldmath$\varepsilon$}}
\newcommand{\bxi}{\bveps}

\newcommand{\bsig}{\mbox{\boldmath$\sigma$}}

\newcommand{\bzeta}{\mbox{\boldmath$\zeta$}}

\def\bz{\mathbf{z}}
 \def\bu{\mathbf{x}}

\def\GG{G}

\def\UU{G}

\def\rank{\mbox{rank}}

\def\bR{{\bf R }}

\def\bL{{\bf L }}
\def\bH{{\bf H}}

\def\alp{\alpha}

\def\be{{\bf e}}
\def\bb{{\bf c}}
\def\half{\frac{1}{2}}

\def\Diag{\mbox{Diag}}

\def\bG{\mathbf{G}}

\def\calS{{\cal S}}

\def\calP{{\cal P}}

\def\barbrho{\bar{\brho}}

\def\eb{\begin{equation}}
\def\ee{\end{equation}}

\newtheorem{thm}{Theorem}
\def\PP{P}

\def\bg{\mathbf{g}}
\def\bx{\mathbf{x}}

\def\calG{{\cal G}}

\def\sig{\sigma}

\def\calX{{\cal X}}

\def\cc{c}
\def\Lam{\Lambda}
\def\calZ{{\cal Z}}
\def\calU{{\cal X}}
\def\be{{\bf e}}
\def\half{\frac{1}{2}}

\def\bv{\mathbf{v}}

\def\vv{v}
\def\bc{\mathbf{c}}
\def\bD{\mathbf{D}}
\def\bK{\mathbf{K}}

\def\Diag{\mbox{Diag}}

\def\bQ{\mathbf{Q}}

\def\bG{\mathbf{G}}

\def\bff{\mathbf{f}}
\def\calS{{\cal S}}
\def\calE{{\cal E}}

\def\calP{{\cal P}}

\def\barbrho{\bar{\brho}}

\def\barbzeta{\bar{\bzeta}}

\def\eb{\begin{equation}}
\def\ee{\end{equation}}
 \newtheorem{remark}{Remark}
 \newtheorem{Corollary}{Corollary}

\begin{document}
\title{Canonical Duality Theory  and Algorithm for  Solving Bilevel  Knapsack Problems with Applications }
\author{{David Yang Gao }
\thanks{Manuscript received ---------; revised September ------------.}%
\thanks{ David Yang Gao is with the Faculty of Science and Technology, Federation University Australia, Mt Helen, Victoria, 3353, Australia (d.gao@federation.edu.au) }}%

\maketitle

\begin{abstract}
A novel canonical  duality theory (CDT)  is presented  for solving general bilevel mixed integer nonlinear optimization  governed by linear and quadratic  knapsack problems.
It shows that  the challenging knapsack problems can be   solved analytically in term of their  canonical dual solutions.  The   existence and uniqueness of these  analytical solutions are  proved.
NP-Hardness of the knapsack problems is discussed.  A powerful CDT algorithm  combined with an alternative iteration  and a volume reduction method
  is proposed for solving the NP-hard  bilevel knapsack  problems. Application is illustrated by a benchmark problem in optimal topology design.
  The performance and novelty of the proposed method are compared with   the popular   commercial codes.
\end{abstract}

\section{Problems and Motivation}
Bilevel optimization is a special kind of optimization where an (upper-level) optimization problem contains another (lower-level) optimization
problem as a constraint.
The first bilevel optimization problem  was  realized in the field of game theory by a German economist Heinrich Freiherr von Stackelberg in 1934 \cite{stack}.
Since then, these problems  have been found  in many real-world applications,
  including  artificial intelligence,   cryptography,  decision science, defensive force structure design,   game theory,     topology  design  of structures and  communication networks,   etc (see \cite{bc,cml,kjy,wz,yhl}).
It was discovered recently by Gao \cite{gao-to17}  that the well-studied  topology  optimization  in computational mechanics should be  formulated as
a   bilevel knapsack problem (BKP), i.e. the upper-level optimization is a knapsack problem, while the lower-level optimization is governed by the well-known minimum total potential energy principle.
The knapsack problem is one of the most fundamental problems in combinatorial optimization. It has been
 studied extensively  for more than a century  in multidisciplinary  fields (see \cite{k-p-p}).
The quadratic knapsack problem (QKP)
  is an extension of knapsack problem that allows for quadratic terms in the
target  function. QKP was  first introduced in 19th century (see  \cite{wit}),
 which has a wide range of applications including telecommunication, transportation network, computer science and economics. In fact, Witzgall first discussed QKP when selecting sites for satellite stations in order to maximize the global traffic with respect to a budget constraint. Similar model applies to problems like considering the location of airports, railway stations, or freight handling terminals \cite{rhy}.
Applications of QKP in the field of computer science is more common after the early days: compiler design problem \cite{h-r-w}, clique problem \cite{p-l-p},    very large scale integration (VLSI) design and manufacturing  \cite{f-m-s}.
In bilevel optimization of  multi-scale complex systems,
both linear and quadratic   knapsack problems can appear in either upper or lower optimization.
Indeed, over the last 20 years,  a variety of bilevel knapsack problems has been proposed (see \cite{bhm,cclw}).

 Generally speaking, bilevel knapsack    problems are  extremely difficult from
the computational point of view and cannot be expressed in terms of classical integer
programs (which can only handle a single level of optimization). Due to the integer constraint and bilevel coupling, bilevel knapsack  problems are  non-differentiable and   non-convex. It has been proven in \cite{vsj}
 that merely evaluating a solution for optimality is also a NP-hard task.
In fact,  it was shown in \cite{hanJS}  that even the linear bilevel problems are  strongly NP-hard.  The proof for the non-existence of a polynomial
time algorithm for linear bilevel problems can be found in \cite{deng}.
Classical  theories and methods  can't be used directly for solving  this type of problems.
Therefore, the general bilevel knapsack problem
 could be one of the most challenging problems in computer science and global optimization.

During the past 40 years  many approximation methods have been proposed for numerically  solving bilevel
optimization  problems. Impressive collections of these methods and applications can be found in the books
 \cite{bard,dem1,z-l-g} as well as the review article \cite{dem2}.
Evolutionary method \cite{smfb,wjl} is a   popular approach, where the lower level optimization
problem is solved corresponding to each and every upper
level member. One of the first evolutionary algorithms for solving bilevel
optimization problems was proposed in the early 1990s.
Mathieu et al. \cite{mpa}  used a nested approach with genetic
algorithm at the upper level, and linear programming at the
lower level. One of the early works on discrete bilevel optimization was
by Vicente et al. (1996) \cite{vsj}, which focused on discrete linear
bilevel programs, and analyzed the properties and existence
of the optimal solution for different kinds of discretization
arising from the upper and lower level variables. The authors
have shown in the paper that certain compactness conditions
guarantee the existence of optimal solution in continuous
linear bilevel programs, discrete-continuous linear
bilevel programs and discrete-discrete linear bilevel programs.
The conditions are equivalent to stating that the inducible
region is non-empty. However, the existence conditions in
the case of continuous-discrete linear bilevel programs are
not straightforward.
In the classical literature, branch-and-bound and branch-and-cut are some of the commonly used deterministic methods to handle discreteness in variables.
But  these well-studied  methods can be computational expensive, and can be used only for solving very small-sized problems. Therefore, global
optimization problems with 200 variables are referred to as ``medium scale", problems with 1000 variables as
``large scale", and the so-called ``extra-large scale" is only around 4000 variables. However, any simple problem
in topology  design can easily have  millions of variables \cite{gao-to18}.
Although  there are some fully polynomial time approximation schemes for discrete optimization,
they are not always applicable in practice due to memory requirements \cite{whs}.
Also these algorithms are based on  heuristic techniques,
 the obtained solutions come with no guarantee of global optimality and may get stuck in local minima.
Therefore,  it  is truly important to develop a powerful deterministic method for solving large-scale  general bilevel  knapsack problems.

Canonical duality theory  (CDT) is a  precise  methodological theory, which can be used not only for modeling complex systems within a unified framework, but also for solving a large class of challenging problems in nonconvex analysis and global optimization \cite{gao-book00}.
This theory is particularly powerful for solving integer constrained problems.
 In 2007  Gao discovered  that  by simply using a standard  canonical measure $\Lam(\bz) = \bz^2 - \bz$, the 0-1 constrained problems can be equivalently converted to a unified
continuous concave maximization problem in a convex  feasible space, which can be solved deterministically via well-developed convex optimization techniques \cite{gao-jimo07}.  This method has been generalized for solving general discrete optimization problems \cite{gao-cace,gao-watsonetal,ruan-gao-minl,wang-etal}.
 Applications of the canonical duality theory to multidisciplinary study
 was given recently \cite{g-l-r-17} and an analytic solution to linear knapsack problem has been obtained \cite{gao-to17,gao-to18}.

This paper deals with a general  bilevel knapsack problem (BKP), in which,  the upper-level optimization is a quadratic knapsack problem coupled with the continuous follower variables, while the lower-level optimization is a mixed integer nonlinear minimization problem involves both leader and follower variables
(since an alternative iteration method is used in this paper, this order of upper and lower levels  can reversed).
The main contributions of this paper include:
\begin{verse}
1. Canonical dual solutions to quadratic   knapsack problems;\\
2. Improved analytic solutions to  knapsack problems and  criteria for NP-hardness;\\
3. A volume reduction method combined with an alternative iteration for solving general BKP;\\
4. A powerful  CDT algorithm for solving general BKP  with an application to topology design.
\end{verse}
The   rest of the paper is organized as follows.
 The next section provides the mathematical formulation of the  bi-level knapsack problem and its special linear case.
 A decoupled alternative iteration is suggested.
 The   theoretical results of this paper are presented in  Section 3, including a canonical penalty-duality method,
 a general  analytic  solution form, as well as the existence and uniqueness of this solution to the quadratic knapsack problem.
In Section \ref{sec4} , the NP-hardness of the knapsack problems is addressed.
Improved analytic solutions are provided for both quadratic and linear knapsack problems.
 Section 5 proposes a volume reduction method and a CPD algorithm for solving  BKP.
 In Section 6, applications are  illustrated by  both 2D and 3D benchmark problems
 in   topology design. The performance and novelty of the CPD algorithm are verified by computational
results. Finally, Section 7 presents conclusions and open problems  for future research.

\section{Bilevel Optimization Modeling and  Knapsack Problems}
The bilevel optimization  problem proposed to solve is given below:
\begin{eqnarray}
 &(\calP):    &   \min_{  \bu \in \calU_a ,   \brho \in \calZ_a} \left\{\PP(\bx,  \bz ) =   \half  \brho^T \bQ( \bu) \bz - \bz^T \bb(\bx)  \right \} \;\;\;\; \label{eq-upp}\\
&& s.t.    \;\; \bu \in \arg \min_{  \bchi \in \calU_a}  \left\{ \Pi( \bchi,  \bz)  = \UU(\bD \bchi, \bz) - \bchi^T \bff  \right \} \;\;\;\;\;\; \;\;\;\; \label{eq-lowp}
\end{eqnarray}
  where $\PP: \calU_a\times \calZ_a \subset \real^m\times \real^n \rightarrow \real $ is  the upper-level target
 function,
  $\Pi :\calX_a \times \calZ_a   \rightarrow \real $ is  the lower-level target
 function, $\bQ(\bx) $ and $\bb(\bx)$ are two pre-given matrix-valued and vector-valued  functions of $\bu \in \calU_a$, respectively; $\bff \in \real^m$ is a given vector;
  $\bD:\calX_a \rightarrow \calG_a \subset  \real^p$ is a linear operator, and $\UU:\calG_a\times \calZ_a \rightarrow \real$ is a given  function, which is usually a nonconvex objective function
   of $\bg = \bD \bx$, but linearly depends on $\bz \in \calZ_a$.
  The upper-level (or leader) variable $\bz$ is a discrete vector, whose domain is a subset of Boolean space (i.e. a knapsack):
 \eb
\calZ_a = \{ \bz = \{ z_i\} \in  \{ 0, 1 \}^n| \;\; \bv^T \bz \le V_c \} ,
\ee
where $V_c > 0$ is a given  size of the knapsack, $\bv = \{ v_i\}  \in \real^n_+$ is a given positive  vector. 
The lower-level (or the follower) variable $\bx$ is a continuous vector, whose domain $\calX_a$ is a convex set of $\real^m$, in which, only linear equality constraints are included.

  Clearly,  the problem $(\calP)$  is a  bilevel knapsack problem.
   Due to the integer constraint, the nonlinearity/nonconvexity of $\Pi(\bx, \bz)$, and the strong coupling between the upper and  lower level  problems,
    the proposed
    problem $(\calP)$ could be one of the most  challenging problems in global optimization and computer science.
    It is fundamentally difficult to solve this problem directly.
    Therefore,  a   decoupled alternative iteration (DAI) method  will be adopted in this paper, i.e., the  coupled bilevel optimization
    is   split into two decoupled subproblems by alternative iteration:
      \begin{verse}
  (i) For a given initial   $\brho_{k-1} \in \calZ_a$, to solve the lower-level problem (\ref{eq-lowp})  first for
  \eb
  \bu_k =  \arg \min \{ \Pi( \bu, \brho_{k-1})   \;\; | \;\; \bu\in \calU_a \} . \label{eq-lowk}
  \ee

  (ii) For the fixed  $\bu_k \in \calU_a$, to solve the upper-level knapsack problem (\ref{eq-upp}) for
  \eb
\bz_k = \arg   \min \{  \PP(   \bu_k, \bz)    | \; \;\; \bz\in \calZ_a \} \label{eq-uppk}
  \ee
  \end{verse}

The upper-level minimization \eqref{eq-uppk}  is the well-known quadratic knapsack problem (QKP).
The  binary variable $\bz = \{ z_i\}$  represents whether item $i$ is included in the knapsack $\calZ_a$.
The  given vector  $\bc(\bx) = \{c_i(\bx) \} \in \real^n_+$ must be
 positive  for any given $\bx \in \calX_a$ such that $c_i$
  is the profit earned by selecting item $i$.
 The given matrix $\bQ(\bx) \in \real^{n\times n} $ must be symmetrical,  its  diagonal elements are usually
$Q_{ii} = 0, \;\; i=1, \dots, n$  and $- Q_{ij} $ is the profit achieved if both item $i$ and $j$ are added.
  An important special case is  the following linear knapsack problem (LKP):
 \eb
 (\calP_{\ell}): \;\; \min \{ \PP_{\ell} (\bz) = - \bb^T \bz | \;\; \bv^T \bz \le V_c,  \;\; \bz \in \{ 0, 1\}^n\}.
 \ee
  A 1998 study of the Stony Brook University Algorithm Repository showed that, out of 75 algorithmic problems, the knapsack problem was the 19th most popular and the 3rd most needed after suffix trees and the bin packing problem.
 It is well-known that the    knapsack problem is NP-hard and there is no known algorithm that can solve the problem in polynomial time.
Actually, even the  linear knapsack problem is listed as one of Karp's 21 NP-complete problems~\cite{karp}.
  In this paper, the  canonical  duality method will be addressed for solving this challenging problem.

For a given design variable $\bz_{k-1}$, the lower-level minimization  \eqref{eq-lowk}
 is the  general nonlinear/nonconvex  optimization problem proposed by the author for multi-scale systems \cite{gao-aip,gao-amma18}.

\begin{remark}[Objectivity,  Modeling, and Canonical Duality]
Objectivity is a central concept in our daily life, related to
reality and truth. In science, the objectivity is often attributed to the property of scientific
measurement, as the accuracy of a measurement can be tested independent from the individual scientist who first
reports it. In continuum physics, a real-valued function $G:\calG_a \rightarrow \real$ is called objective if and only if it is an invariant under orthonormal transformation (see Definition 6.1.2 \cite{gao-book00}), i.e.
\eb
G(\bR \bg) = G(\bg) \;\; \forall \bg \in \calG_a , \;\; \forall   \bR \in SO(p) , 
\ee
where $SO(p) = \{ \bR \in \real^{p\times p}  | \;\;  \bR^T = \bR^T, \; \det \bR =1 \}$ is a special orthogonal  group in $\real^p$.
Physically speaking, an objective function is governed by the intrinsic physical law of the system,
which does not depend on observers. Therefore, the objectivity is essential for any real-world
mathematical models. In continuum physics, it  is also called the principle of material
frame indifference \cite{oden}.  Geometrically, an objective function does not depend on rigid rotation of the
system considered, but only on certain measure of its variable. In the Euclidean space $\calG_a \subset \real^p$, the simplest objective function is the $\ell_2$-norm
$\| \bg \|$ since $\| \bR \bg\|^2 = \bg^T \bR^T \bR \bg = \| \bg \|^2 \;\; \forall \bR \in SO(p)$. Therefore, an objective function must be nonlinear.

Correspondingly, the linear term $F(\bx) = \bx^T \bff  $ in the lower-level problem \eqref{eq-lowp}  is called a {\em subjective function} \cite{gao-aip,gao-to18}, where the given input $\bff$ and the constraints in $\calX_a$   depend only   on each given problem.
Thus, if the operator ``$-$" is considered as the predicate, the
 difference  $\Pi(\bx) = \UU(D \bx) - F(\bx)$ between  object and subject
 forms a complete and grammatically correct mathematical formulation. 
  In continuum physics, if $\GG(\bg)$ is the free (or internal) energy, $F(\bx)$ is the input (or external)
 energy, then $\Pi(\bx)$ is the total potential energy and the minimum total potential energy principle leads to
a general variational/optimization problem in mathematical physics \cite{gao-strang89}.
By the fact that $\calG_a \subset \real^p$ and $\calX_a \subset \real^m$ can be  in different dimensional spaces with different measures,
  the lower-level minimization
$\min \Pi(\bx) $ presented in this paper covers general constrained nonconvex/nonsmooth  optimization
problems in multi-scale complex systems \cite{gao-book00,gao-amma18,g-l-r-17}.

According to  \cite{ciarlet},  $\UU(\bg)$ is an objective function if and only if
there exists an objective measure $\beps = \bg^T \bg$ and a real-valued function $\Phi(\beps)$ such that $G(\bg) =\Phi(\bg^T \bg)$.
In continuum physics  and differential geometry, if $\bg = D \bx$ is the deformation gradient, then the objective measure $\beps = \bg^T \bg$ is the well-known right Cauchy-Green tensor. By the fact that the  free energy $\Phi(\beps)$ is usually convex (say the St Venant-Kirchhoff material \cite{gao-book00}),
 the duality relation $\beps^* = \nabla \Phi(\beps)$ is invertible. This one-to-one constitutive relation is called the canonical duality and $\Phi(\beps)$ is called the canonical function. These basic truths in continuum physics laid a foundation for the
canonical duality theory, in which, the key idea  of the
{\em canonical transformation} is to choose a nonlinear operator (not necessary to be objective)
$\beps= \Lam(\bx)$ and a canonical function $\Phi(\beps)$ such that
the  given nonconvex function $G(D\bx)$ can be written in the canonical form $G(D\bx) = \Phi(\Lam(\bx))$.
By the fact that the  objectivity plays a fundamental role in mathematical modeling,   the canonical
duality theory can be powerfully used for solving
many  challenging problems in multidisciplinary fields \cite{g-l-r-17}.

 However, this  important concept of the objectivity  has been extensively  misused in optimization literature such that the general problem in nonlinear optimization (or programming)
has been  formulated as
\eb
\min f(x) \;\; s.t. \;\;  g(x) = 0, \;\; h(x) \le 0 ,
\ee
where $f(x)$ is called the ``objective function"\footnote{This terminology is used mainly in English literature. The function $f(x)$
 is correctly called the target function in all Chinese and Japanese literature, or the goal
function in some Russian and German literature by serious researchers.}, which is allowed to be any arbitrarily given function, even a linear
function. Clearly, this mathematical problem is too abstract. Although it enables one to model a very wide range of
problems, it comes at a price: many global optimization problems are considered to be NP-hard. Without detailed
information on these arbitrarily given functions, it is impossible to have a powerful theory for solving the artificially
given constrained problem. Also, due to this conceptual mistake,   the canonical
duality theory has been mistakenly challenged 
in several publications\footnote{See  M.D. Voisei, C. Zalinescu, Some remarks concerning Gao-Strang’s complementary gap function, {\em Appl. Anal.} 90 (2011) 1111-1121.}.
The conceptual
mistake  in these false challenges revealed a big gap   between physics and optimization.
Interested readers are recommended to read the recent papers \cite{gao-opl16,g-r-l-17} for details.  \hfill $\diamondsuit$ 
\end{remark}

 The canonical duality theory for solving the lower-level  nonconvex continuous optimization problem (\ref{eq-lowk})
has been studied extensively during the past 15 years    (see \cite{gao-jogo00, gao-cace,g-l-r-17}).
 As long as the nonconvex function $\UU(\bD \bx, \bz)$ in $\Pi(\bx, \bz)$  can be written in the
canonical form $\Phi(\Lam(\bx), \bz)$, this problem can be  solved easily by the canonical duality theory
  to obtain both global and local optimal solutions (see many real-world applications in \cite{g-l-r-17}).
This paper will show how to use the canonical duality theory for solving the challenging upper-level knapsack problem \eqref{eq-uppk}.

 \section{Canonical  Duality  Solution to Upper-Level Knapsack Problem}
\label{Pure}
 For a fixed $\bx \in \calX_a$,  the upper-level  optimization (\ref{eq-uppk})  can be written in  the typical  quadratic knapsack  problem:
  \eb
  (\calP_q): \;\; \min_{  \brho\in \{ 0,1\}^n } \left\{ \PP_q(\brho) = \half    \brho^T \bQ \bz  - \bb^T \brho  \;\; |\;\; \bv^T \brho \le V_c  \right\}.
  \ee

    According to the canonical duality theory for mathematical modeling \cite{gao-opl16,gao-to17}, the inequality
    $\brho \cdot \bv \le V_c$ is a geometrical constraint, while the  integer constraint $\brho \in \{ 0,1\}^n$ in   $(\calP_q)$
 is a constitutive   condition \cite{gao-to18,gao-amma18}, which can be equivalently replaced by  the canonical constraint
 $\bz^2 - \bz = {\bf 0} \in \real^n \;\;\forall \bz \in \real^n$, where $\bz^2 = \bz\circ \bz  = \{ z_i^2\}  \in \real^n$.
 Therefore, by introducing a penalty parameter $\beta > 0$ and let  $\calZ_b = \{    \bz \in \real^n  \; | \;\; \bv^T \bz \le V_c  \}$,  the problem $(\calP_q)$ can be relaxed by  the canonical  penalty function:
\eb
(\calP_\beta): \;\;
\min_{ \bz \in \calZ_b} \left\{ \PP_\beta(\bz) = \half \bz^T \bQ \bz +  \half \beta \|\bz^2 - \bz \|^2  -  \bb^T  \brho 
\right\} \label{eq-pda} .
\ee
Clearly, we have
\eb
\{ \PP_q(\bz)  | \;\; \bz \in \calZ_a\} = \lim_{\beta \rightarrow \infty} \{   \PP_\beta(\brho )   | \;  \bz \in \calZ_b  \}.
\ee
Although the integer constraint is relaxed by the   penalty function,
 \eb
 \Psi_\beta(\brho) =  \half \beta \|\bz^2 - \bz \|^2   ,
 \ee
  the    problem (\ref{eq-pda}) is a nonconvex minimization in $\real^n$, which
 is still considered to be NP-hard   by traditional methods. This is the reason why the traditional external penalty method can be used mainly for linear constrained problems \cite{l-g-opl}.
  To solve this problem by  the canonical duality theory,  we need first to introduce the following canonical transformation:
  \[
  \Psi_\beta(\brho) = \Phi_\beta(\Lam(\brho)),     
\]
\[
   \Phi_\beta(\bxi) =   \half  \beta \| \bxi\|^2, \;\;
 \bxi = \Lam(\brho) = \brho^2 - \bz  \in    \real^n.
  \]
Clearly,  $\Phi_\beta: \real^n  \rightarrow \real$ is a convex quadratic function. Its canonical dual can be simply given by the Legendre transformation
\eb
\Phi^*_\beta(\bsig) = \max_{\bxi \in \calE } \{ \bxi^T \bsig - \Phi_\beta(\bxi)   \} =
\half \beta^{-1}   \|\bsig\|^2 .
\ee
 Thus, replacing $\Psi_\beta(\brho) $ in $P_\beta(\brho )$  by the Fenchel-Young equality $\Phi_\beta(\Lam(\brho) = \Lam(\brho)^T \bsig - \Phi^*_\beta(\bsig)$,
 and introducing a Lagrange multiplier $\tau \ge 0$,
the canonical penalty function $\PP_\beta(\bz)$ can be reformed as  the   Gao-Strang total complementary function
  \eb
   \Xi_\beta(\brho, \bsig,\tau) =  \half   \brho^T \bG(\bsig ) \bz  -  \half  \beta^{-1} \| \bsig\|^2
  -  \brho^T \btau (\bsig, \tau  )  -  \tau V_c ,
\ee
  where
\eb
\bG (\bsig) = \bQ + 2 \Diag(\bsig), \;\; \btau (\bsig, \tau  ) = \bb  - \tau \bv   + \bsig.
\ee
 By introducing a  canonical dual feasible  space
\eb
 \calS_a^+ = \{ \bzeta = ( \bsig, \tau) \in \real^{n+1} | \;\; \bG(\bsig)  \succ 0 , \;\; \tau >  0 \},
\ee
 the function  $  \Xi_\beta(\brho, \bzeta) $  is convex in $\brho \in \real^n$ for   any given $\bzeta \in \calS^+_a$.
The canonical penalty-duality function $ \PP^d_\beta : \calS^+_a \rightarrow \real$  can be defined by
\begin{eqnarray*}
\PP_\beta^d(\bzeta) &=& \min_{\brho \in \real^n} \{  \Xi_\beta (\brho, \bzeta) |\;\; \bzeta \in \calS^+_a \} \\
&=&  - \half  \btau(\bzeta)^T \bG(\bsig) ^{-1} \btau(\bzeta) -
\half \beta^{-1}  \| \bsig  \|^2 -\tau V_c .
\end{eqnarray*}
Thus, the  canonical penalty-duality  problem  can be proposed in the following
 \eb
 (\calP^d_\beta): \;\;\;\; \max \left\{ \PP^d_\beta (\bzeta) | \;\;\bzeta \in \calS^+_a \right\}.
 \ee
  \begin{thm}[Complementary-Dual Principle]\label{thm1}
 For any given $\beta  > 0$, if $  (\bz_\beta,  \bzeta_\beta)    \in  \real^n \times   \calS^+_a$ is a KKT point of $\Xi_\beta(\bz, \bzeta)$,     then
 $\bz_\beta$  is a global minimum solution to  the canonical  penalty   problem $(\calP_\beta)$,
 $\bzeta_\beta$ is a  solution to $(\calP^d_\beta)$, and
  \eb\label{eq-cdp}
P_\beta(\bz_\beta) = \min_{\bz \in \calZ_b} \PP_\beta(\bz) =    \Xi_\beta (\barbrho, \barbzeta )=
\max_{\bzeta \in \calS^+_a} \PP^d_\beta(\bzeta) =    P_\beta^d(\bzeta_\beta).
\ee
\end{thm}
 {\bf Proof}. By the fact that $\Xi_\beta :\real^n \times \calS^+_a \rightarrow \real$ is a saddle function,  if $  (\bz_\beta,  \bzeta_\beta)    \in  \real^n \times   \calS^+_a$ is a KKT point of $\Xi_\beta(\bz, \bzeta)$,   it must be a saddle point of $\Xi_\beta(\bz, \bzeta)$.
 Then by the definition of $\PP_\beta^d(\bzeta)$, we have
 \begin{eqnarray} 
\Xi_\beta(\bz_\beta, \bzeta_\beta) &=&  \min_{\bz \in \real^n} \max_{\bzeta \in \calS^+_a} \Xi(\bz, \bzeta) =   \max_{\bzeta \in \calS^+_a}  \min_{\bz \in \real^n}\Xi(\bz, \bzeta) \nonumber  \\
&  =&  \max_{\bzeta \in \calS^+_a}  \PP^d_\beta(\bzeta) = \PP^d_\beta(\bzeta_\beta).
 \end{eqnarray}
   By the KKT conditions 
   \[
\bv^T \bz_\beta  \le V_c, \;\; \tau_\beta \ge 0 , \;\;
 \tau_\beta(  \bv^T \bz_\beta -  V_c) = 0 ,
 \]
and the condition $\bzeta_\beta = (\bsig_\beta, \tau_\beta) \in \calS^+_a$,
we have  $\tau_\beta > 0$ and  $\bv^T \bz_\beta -  V_c = 0$.
By the
   convexity of the penalty function $\Psi_\beta (\beps)$, we have
$ \Psi^{**}_\beta (\beps) = \Psi_\beta(\beps) $.
 Thus, for any given $\beta > 0$
\eb
\PP_\beta(\bz_\beta) = \min_{\bz \in \calZ_b} \PP_\beta (\brho ) = \min_{\bz \in\calZ_b} \max_{\bzeta \in \calS^+_a}  \Xi_\beta (\brho, \bzeta)   = \Xi_\beta (\bz_\beta, \bzeta_\beta) .
\ee
Thus, $\bz_\beta$ is a global optimal solution to $(\calP_\beta)$.  \hfill $\Box$\\

Theorem 1 shows that the nonconvex minimization problem $(\calP_\beta)$ is canonically dual to a concave maximization problem  $(\calP^d_\beta)$ in
a convex space $\calS^+_a$, which can be  solved easily by well-developed convex minimization methods.

 \begin{remark}[Canonical Penalty-Duality    and $\beta$-Perturbation]
By the facts that
\eb\label{eq-pane}
  \Phi(\bxi) = \lim_{\beta \rightarrow \infty}  \Phi_\beta(\bxi)    = \left\{ \begin{array}{ll}
0 & \mbox{ if } \bxi  = {\bf 0} \in \real^n\\
+ \infty & \mbox{ otherwise} \end{array} \right.
 \ee
and
\eb\label{eq-pane*}
\Phi^*(\bsig) = \sup_{\bxi \in \real^n}  \{ \bxi^T \bsig - \Phi(\bxi) \} =
0  \;\; \forall \bsig \in \real^n,
 \ee
we have
\begin{eqnarray*}
\PP^d_q(\bzeta)  & = & \lim_{\beta \rightarrow + \infty} \PP_\beta^d (\bzeta) \\
&=&  - \half  \btau(\bzeta)^T \bG(\bsig)^{-1} \btau(\bzeta)
  -\tau V_c  - \Phi^*(\bsig),
  \end{eqnarray*} 
  which is exactly the canonical dual function to the primal function $\PP_q(\bz)$.
  Therefore, $(\calP^d_\beta)$ is actually the so-called $\beta$-perturbation of the canonical dual problem $(\calP^d_q)$:
  \eb
  (\calP^d_q): \;\; \max \{ \PP^d_q (\bzeta) | \;\; \bzeta \in \calS^+_a \} .
  \ee
  The penalty-duality method  was first proposed  by  Gao for solving convex variational problems \cite{gao-thesis,gao-cs88}.
 The $\beta$-perturbation method for nonconvex  integer constrained problems was  first proposed
 in \cite{gao-to17,gao-ruan-jogo10}.
 It was proved by   Theorem 7  in \cite{gao-ruan-jogo10} that there exists a $\beta_c>0$ such that for any given $\beta \ge \beta_c$,
both  the $\beta$-perturbed canonical dual problem
$(\calP^d_\beta)$
 and the problem $(\calP^d_q)$ have  the same solution set.
 This shows the relation between the canonical penalty-duality method proposed in this paper and the $\beta$-perturbation method proposed in \cite{gao-ruan-jogo10}. \hfill$\diamondsuit$ 
 \end{remark}

  \begin{thm}[Analytic  Solution to Knapsack Problem]\label{thm2}
 For any given  $\beta  > 0$, if  $   \bzeta_\beta =(\bsig_\beta, \tau_\beta)  \in   \calS^+_a$ is a   solution to $(\calP^d_\beta)$,
   then  the  vector defined by
 \eb \label{eq-solu}
 \bz_\beta   =   \bG(\bsig_\beta)^{-1} \btau (\bzeta_\beta)
 \ee
is a global minimum solution to the canonical penalty   problem $(\calP_\beta)$.

Moreover, there exists   $\beta_c  \gg 0$ such that $\beta\ge \beta_c$ and $\bz_\beta \in \calZ_a$, then $\bz_\beta $ is a global optimal solution to the knapsack problem $(\calP_q)$
and
 \eb\label{eq-cdp}
P_q(\bz_\beta) = \min_{\bz \in \calZ_b} \PP_\beta(\bz) =    \Xi_\beta (\barbrho, \barbzeta )=
\max_{\bzeta \in \calS^+_a} \PP^d_\beta(\bzeta)  =   P_q^d(\bzeta_\beta).
\ee
\end{thm}
 {\bf Proof}. It is easy to verify  that if  $(\bz_\beta, \bzeta_\beta)$ is  a KKT point of $\Xi_\beta(\bz, \bzeta)$, the
 criticality condition $\nabla_\bz \Xi_\beta(\bz_\beta, \bzeta_\beta) = 0 $ leads to  $ \bz_\beta   =   \bG(\bsig_\beta)^{-1} \btau (\bzeta_\beta) $.
 By Theorem \ref{thm1} we know that if  $   \bzeta_\beta =(\bsig_\beta, \tau_\beta)  \in   \calS^+_a$ is a   solution to $(\calP^d_\beta)$, the vector defined by (\ref{eq-solu})
 must be a solution to $(\calP_\beta)$.
Since the penalty function $\Phi_\beta(\bz) \ge 0 \;\; \forall \bz \in \real^n$, there must exists a sufficiently big $\beta_c \gg  0$
such that $\Phi_\beta(\bz_\beta)  = 0 \; \forall \beta \ge \beta_c$, i.e.   $\bz_\beta \in \calZ_a$.
Thus, $\PP_q(\bzeta_\beta) = \PP_\beta(\bz_\beta) = \min \{ \PP_q(\bz) | \; \bz \in \calZ_a\}$. \hfill $\Box$\\

 Since $\calS^+_a$ is a convex open set, its boundary is
 \eb
 \partial \calS^+_a = \{ \bzeta = (\bsig, \tau)  \in \real^{n+1} | \;\;  \det  \bG(\bsig) = 0 , \;\; \tau = 0 \}.
 \ee
 The following theorem is important for understanding the NP-hardness of the knapsack problem.
  \begin{thm}[Existence and Uniqueness for $(\calP_q)$]\label{thm3}
 For any given $\bb  \in \real^n$, $ \bv \in \real^n_+ $, and $\bQ \in \real^{n\times n}$,   if $\calS^+_a$ has a  non-empty interior $\bzeta_o \in \calS^+_a$ such that
 $\PP_q^d(\bzeta_o) >  - \infty $ and
 \eb
 \lim_{\bzeta \rightarrow \partial \calS^+_a } \PP^d_q(\bzeta) = - \infty \;\; \;\; \forall \bzeta \in \calS^+_a ,\label{eq-corc}
 \ee
 then  the  problem $(\calP_q^d)$ has a unique solution $\bzeta_c \in \calS^+_a$ and $\bz_c= \bG(\bsig_c) ^{-1}\btau(\bzeta_c) $ is a unique global optimal solution to  $(\calP_q)$.
 \end{thm}
 {\bf Proof}. This theorem is a direct application of the canonical duality  theory (see Theorem 4 in  \cite{gao-cace}\footnote{The condition $\PP^d(\bsig_o) < -\infty$ in  \cite{gao-cace}  is a typo, it should be  $\PP^d(\bsig_o) > -\infty$.}).
  By the facts that $\PP^d_q(\bzeta)$ is concave and its feasible set $\calS^+_a$ is convex and nonempty,   if the condition (\ref{eq-corc})
 holds, then   the canonical penalty-dual function $\PP^d_q:\calS^+_a \rightarrow \real $ is strictly concave, bounded up and
$-\PP^d_q(\bzeta)$ is coercive on $\calS^+_a$. Therefore, $(\calP^d_q)$ has a unique solution $\bzeta_c$. By Theorem \ref{thm2}, we know that the associated  $\bz_c= \bG(\bsig_c) ^{-1}\btau(\bzeta_c) $ is a unique global optimal solution to  $(\calP_q)$. \hfill $\Box$ \\

 As an important application, let us re-call  the   linear knapsack problem:
 \eb
 (\calP_{\ell}): \;\; \min_{\bz \in \calZ_a}  \{ \PP_{\ell} (\bz) = - \bb^T \bz  \} .
 \ee
The associated canonical penalty problem is
\eb
(\calP_{\ell\beta}): \;\;  \min_{  \bz \in \calZ_b}  \left\{ \PP_{\ell\beta} (\bz) = \half \beta \|\bz^2 - \bz \|^2   - \bb^T \bz \;\;
  \right\} .
\ee
 Since the matrix $\bQ $ vanished in this case, the  canonical penalty-dual function  has the following simple form:
 \eb
  \PP^d_{\ell\beta} (\bzeta) = -  \half  \sum_{i=1}^n   \left( \half  \sig^{-1}_i (\sig_i + c_i - \tau v_i)^2 + \beta^{-1} \sig_i^2 \right)   - \tau V_c  .
 \ee
Let
$ \bsig^{-1} = {  \Diag}(\bsig)^{-1}  =\{ \sig^{-1}_i\}$ and 
\[
\calS^+_a = \{ \bzeta = (\bsig, \tau) \in \real^{n+ 1} | \;\; \bsig > {\bf 0} , \;\; \tau > 0 \}.
\]
The   canonical penalty-dual problem for the linear knapsack problem is 
\eb
(\calP_{\ell\beta}^d): \;\; \max \{   \PP^d_{\ell\beta} (\bzeta) | \;\; \bzeta \in \calS^+_a\}.
\ee

 \begin{Corollary}[Analytic  Solution to Knapsack Problem \cite{gao-to18}]\label{thm4}
 For any given $\bb, \bv \in \real^n_+$ and $V_c , \beta > 0$,
 if    $   \bzeta_\beta =(\bsig_\beta, \tau_\beta)  \in   \calS^+_a$ is a   solution to $(\calP^d_{\ell\beta})$
 and $\tau_\beta \bv - \bb \neq {\bf 0}$,
   then  the penalty   problem $(\calP_{\ell\beta})$ has a unique global minimum solution which is given analytically
   by
 \eb \label{eq-solul}
 \bz_\beta   =  \half  \bsig_\beta^{-1} \circ   (\bsig_\beta + \bb - \tau_\beta \bv) .
 \ee
Moreover, there exists a $\beta_c \gg  0$ such that $\beta\ge \beta_c$ and  $\bz_\beta \in \calZ_a$, then $\bz_\beta $ is a global optimal solution to the linear knapsack problem $(\calP_\ell)$.
\end{Corollary}

This special  result for linear knapsack problem was obtained recently in \cite{gao-to18}.  For $\beta = 0$, the canonical dual problem $ (\calP_{\ell}^d)$ of $(\calP_\ell)$ is 
 \eb
\max_{\bzeta \in \calS^+_a } \left\{ \PP^d_{\ell } (\bzeta) = -  \half  \sum_{i=1}^n   \left( \half  \sig^{-1}_i (\sig_i + c_i - \tau v_i)^2   \right)   - \tau V_c \right\} .
 \ee
It is easy to prove that if $\tau_\beta \bv - \bb \neq {\bf 0}$, the condition (\ref{eq-corc})  is reduced by
 \eb
 \lim_{\bsig \rightarrow {\bf 0}^+ } \PP^d_{\ell} (\bsig, \tau_\beta) = -\infty  .
 \ee
 Thus, by Theorem \ref{thm3} we know that the vector (\ref{eq-solul}) is a unique solution to  the problem $(\calP_{\ell\beta})$ and it is a solution to the  linear knapsack problem
 $(\calP_\ell)$ if $\beta \ge \beta_c \gg 0$ such that $\bz_\beta \in \calZ_a$.

Actually, for any given $\beta > 0$, the criticality condition $\nabla \PP^d_{\ell\beta} = 0$ leads to the following algebraic equations:
\eb
  \beta^{-1} \sig_i^3 + \sig_i^2 = (\tau v_i - c_i)^2, \;\; i = 1, \dots, n, \label{eq-cdas}
\ee
\eb
\sum_{i=1}^n \half \frac{v_i}{\sig_i} ( \sig_i- v_i \tau + c_i ) - V_c = 0 .\label{eq-cdv}
\ee
 It was  proved in \cite{gao-book00,gao-jimo07} that for any given $\beta > 0$,  $\tau \ge 0$ and $\bb \in \real^n$ such that
 $\theta_i = \tau v_i - c_i \neq 0, \forall  i = 1, \dots, n$, the canonical dual algebraic equation (\ref{eq-cdas}) has a unique
positive real solution
\eb
\sigma_i  =  \frac{1}{6} \beta   [- 1 +  \phi_i(\tau ) + \phi_i^c(\tau )] > 0 , \;\; i = 1, \dots, n
\label{eq-solus}
\ee
where
\[
\phi_i(\vsig )  = \eta^{-1/3} \left [2 \theta_i^2 - \eta + 2 \theta_i  \sqrt{ \theta_i^2 -\eta } \right]^{1/3} ,
 \;\; \eta = \frac{4 \beta^2}{27},
\]
and $\phi_i^c $ is the complex conjugate of $\phi_i $, i.e. $\phi_i  \phi_i^c  = 1$.

On the other hand, for a given $\bsig \in \real^n_+$, the Lagrange multiplier $\tau$ can be uniquely obtained by
\eb
        \tau = \frac{ \sum_{i=1}^n  v_i (1 + c_i  /\sigma_i) - 2 V_c }{\sum_{i=1}^n v_i^2/  \sigma_i}. \label{eq-solust}
\ee
 It is easy to prove that for any given
$\beta > 0$ and $\tau_{k-1} > 0$,
the solution $\bsig_k $ produced by \eqref{eq-solus} satisfies
$\bsig_k < \bb$.
Thus, for a given initial $\tau_0 > 0$, an alternative iteration can be used for solving \eqref{eq-solus} and \eqref{eq-solust} and the
   sequence  $(\bsig_k, \tau_k) \in \calS^+_a$ approaches to the global optimal solution of $(\calP^d_\ell)$ in polynomial time.

Theorem \ref{thm4} shows that although the canonical dual problem is a concave maximization in continuous space,  it produces the analytical solution  (\ref{eq-solu})
to  the well-known integer Knapsack problem $(\calP_q)$. This truth was first discovered  by Gao in 2007
for  general quadratic integer programming problems (see Theorem 3, \cite{gao-jimo07}).

\section{Improved Solutions and NP-Hardness  } \label{sec4}
 By the fact that $\balp \circ \bz^2 = \balp\circ  \bz \;\; \forall \bz \in \{0,1\}^n, \;  \forall \balp \in \real^n$,
for any given symmetrical  $\bQ \in \real^{n\times n}$ we can choose an $\balp$ such that   $\bQ_\alp= \bQ +2  \Diag (\balp) \succeq 0$. Thus, by $\bc_\alp = \bc +  \balp $,
the problem $(\calP_q)$ can be equivalently written in the  so-called $\alp$-perturbation form  \cite{gao-cace}:
\[
(\calP_\alp): \;\;    \min_{ \brho\in \{ 0,1\}^n}
\left\{ \PP_{\alp} (\brho) = \half    \brho^T \bQ_\alp  \bz  - \bb_\alp^T \brho  \;\; |\;\; \bv^T \brho \le V_c   \right\} .
  \]
Since  $\rank \bQ_\alp = r \le n$,   there must exist  (see \cite{strang}) a
 $\bL  \in \real^{r\times n} $ and
$\bH \in \real^{r\times r}$ with $\rank \bL = \rank \bH = r$ and $\bH \succ  0$
 such that $\bQ_\alp = 4 \bL^T \bH\bL$.
Similar to  the $\alp$-perturbed canonical dual problem $(\calP^g_{ip})$ given in \cite{gao-cace},
 the  canonical dual function in $(\calP^d_q)$ can be
reformulated as \cite{gao-amma18}:
\[
  \PP^g_{\alp}   (\bsig,\tau) = - \half  \mbox{Abs}[ \bphi(\bsig, \tau)]
    - \half \bsig^T \bH^{-1} \bsig - \tau V_b  + d  ,
\]
where $V_b = V_c - \half \sum_{i=1}^n v_i , \;\; d = \frac{1}{8}  \sum_{i=1}^n  (2 \alp_i+ \sum_{j=1}^n Q_{ij})  - \half \sum_{i=1}^n (\cc_i  + \alp_i) $ are two constants, and
\eb
   \bphi ( \bzeta ) =     \bc  -   \tau \bv -  2  \bL^T \bsig  - \half  \bQ \be  , \;\; \be = \{ 1\}^n .
\ee
The notation $\mbox{Abs}[ \bphi(\bsig,\tau)]$ denotes
$ \mbox{Abs}[ \bphi(\bsig,\tau) ]= \sum_{i=1}^n | \phi_i(\bsig,\tau) | $.
Let 
\[
\calS^+_c= \{ \bzeta =(\bsig, \tau) \in \real^{r+1} | \;\;  \;\; \tau >  0 \}.
\]
Then the improved canonical dual problem to $(\calP_\alp)$ can be proposed as
\eb
(\calP^g_{\alp}): \;\;\; \max \left\{ \PP^g_{\alp}   (\bsig,\tau) | \;\; \bzeta \in \calS^+_c \right\}.
\ee
In many real-world applications, we have $r\ll n$, thus the problem $(\calP^g_{\alp})$ is much easier than $(\calP^d_q)$. Similar to Theorem 7 in \cite{gao-cace}, we have the following improved result.
 \begin{thm}[Improved Analytic  Solution to $(\calP_q)$]\label{thm-eqkp}
 For any given $V_c > 0$,   $\bv, \bc  \in \real^n_+ $, $\bQ \in \real^{n\times n}$ and an $\balp \in \real^n $ such that
  $\bQ_\alp = \bQ + 2 \Diag (\balp)  = 4 \bL^T \bH\bL$,  and
$\bH \succ  0$, if
 $  \bzeta_c =\{ \bsig_c,   \tau_c\} $ is a  solution to $(\calP^g_{\alp})$ and
\eb
\phi_i( \bzeta_c)   \neq { 0} \;\; \forall i \in \{ 1, \dots, n \} ,\label{eq-exis}
\ee
then the quadratic knapsack problem $(\calP_{q})$ has a unique  global optimal solution
\eb\label{eq-solu2}
\bz_c = \half \left\{ \frac{  \phi_i (  \bzeta_c) }{| \phi_i (  \bzeta_c) |} + 1 \right\}^n
\ee
and
 \eb
 \PP_{q} (\bz_c) = \min_{\brho \in \calZ_a} \PP_{q} (\brho) = \max_{\bzeta \in \calS^+_c}  \PP^g_{\alp} (\bzeta)
=   \PP^g_{\alp}(\bzeta_c).
\ee
  Otherwise,  if $\phi_i( \bzeta_c) = 0 $ for at least one $i \in \{1, \dots, n\}$, then
$(\calP_{q})$  has at least two  solutions.
\end{thm}

It is easy to prove that $\PP^g_\alp(\bzeta)$ is strictly concave on the convex feasible set
$\calS_c^+$. If $ \bzeta_c \in \calS^+_c$ is a solution to $(\calP^g_\alp)$, this solution must be unique.
 Since $\dim  \calS^+_c = r+1 \le n+1$, the problem $(\calP^g_\alp)$ is easier than $(\calP^d_\beta)$.
 In the case that  $(\calP^g_\alp)$ does not have a solution to satisfy  \eqref{eq-exis},
we can oppositely  chose $\balp \in \real^n$  such that $\bH \prec 0$. But in this case, $ \PP^g_\alp(\bzeta)$ is a d.c. function (difference of convex functions) and the corresponding problem should be \cite{gao-cace,gao-watsonetal}
\eb
(\calP^g_\alp): \;\; \min \mbox{sta}  \left\{ \PP^g_\alp(\bzeta) | \;\; \bzeta \in \calS^+_c \right\},
\ee
where $ \min \mbox{sta} \{ f(x)\} $ means to find the minimum stationary point of $f(x)$.
This is a nonsmooth d.c. programming problem.
A so-called  VTDIRECT parallel algorithm has been used successfully by Gao {\em et al}
 \cite{gao-watsonetal} for finding global optimum solutions to quadratic integer programming problems, but not in polynomial time.

For the  linear Knapsack problem, the improved canonical dual has a very simple form:
\eb
(\calP^g_{\ell }): \;\;\; \max_{ \tau \ge 0} \left\{  \PP^g_{\ell }(\tau) = -\half \sum_{i=1}^n ( | \cc_i - \tau v_i |  -  \tau v_i)  - \tau V_c  \right\}. \label{eq-lknap}
\ee
\begin{Corollary}[Improved Analytic  Solution  to $(\calP_\ell)$]
For any given  $V_c > 0$, $\bv, \bc  \in \real^n_+ $,    if $\tau_c >  0$ is a solution to $(\calP^g_{\ell })$
and $\theta_i= \tau_c   v_i - \cc_i  \neq 0 \;\forall i\in \{ 1,\dots, n\}$, then the linear knapsack problem $(\calP_\ell)$
has a  unique  global optimal solution
\eb \label{eq-rtau}
\bz_c= \half \left\{  \frac{\cc_i - \tau_c  v_i }{|\cc_i - \tau_c   v_i |} + 1 \right\}^n
\ee
 and
$ \PP_{\ell} (\bz_c) = \PP^g_{\ell } (\tau_c)$.
\end{Corollary}

Actually, from the proof of Corollary \ref{thm4} we know that if there exists a $\tau_\beta > 0$ such that
 $ \theta_i = \tau_\beta \vv_i - \cc_i\neq 0 \;\; \forall i \in \{ 1, \dots, n\}$, the canonical dual algebraic equation \eqref{eq-cdas} has a unique positive solution $\bsig_\beta$
\cite{gao-book00,gao-na00,gao-jimo07}.
For a sufficiently big $\beta \gg 0$ we have $\bsig_\beta = \bsig_c$ and
 the solution $\tau_\beta$ by \eqref{eq-solust} is exactly the solution $\tau_c$ of the problem
 $(\calP^d_{\ell\alp})$ (see \cite{gao-to18}). Therefore, $(\calP_\ell)$ has a unique solution defined by either \eqref{eq-solu} or \eqref{eq-rtau}.

\begin{remark} [Criteria for Non-Uniqueness and NP-Hardness]
Theoretical results presented so far show that if  the condition \eqref{eq-corc}  holds, the canonical dual problem $(\calP_q^d)$ has a unique solution
$\bzeta_c = (\bsig_c, \tau_c) \in \calS^+_a$, which can be obtained deterministically in polynomial time since  $(\calP_q^d)$ is equivalent to a convex minimization problem.
 In this case, both the quadratic and the  linear  knapsack problems are  not NP-hard and  their solutions $\bz_c$ can be analytically given by Theorems \ref{thm2} and Corollary \ref{thm4}, respectively. Also, we must have $\bz_c^T \bv = V_c$ since $\tau_c > 0$.
 For the linear knapsack problem $(\calP_\ell)$, the  condition \eqref{eq-corc} is simply $\tau_c \bv \neq \bc$.

On the other hand, if the canonical dual problem $(\calP_q^d)$ has no solution in $\calS^+_a$, the primal problem $(\calP_q)$ could be   NP-hard,
which is a conjecture first proposed by Gao in 2007 \cite{gao-jimo07}, i.e.
\begin{verse}
{\bf Conjecture of NP-Hardness}. A global optimization problem is NP-hard only if its canonical dual has no solution in $\calS^+_a$.
\end{verse}
 It is also an open problem left in \cite{gao-cace,gao-ruan-jogo10}.
The reason for   NP-hard problems and possible   solutions  were discussed recently in \cite{gao-aip,gao-to18}.
For the linear knapsack problem $(\calP_\ell)$, as long as  $\tau_c \vv_i = \cc_i$ for any one $i \in \{1, \dots, n\}$,
  the canonical dual algebraic equation \eqref{eq-cdas} has at least  two repeated solutions $\bsig_c$ located on the boundary of $\calS^+_a$.
In this case, the primal  problem $(\calP_\ell)$ has multiple solutions, i.e. it is not well-posed \cite{gao-aip,gao-to18}.
A linear perturbation method for   solving this case was proposed recently in \cite{chen-gao-amma,gao-to18,wang-etal}.
\end{remark}

\section{Volume Reduction Method and Canonical  Duality Algorithm}

Theoretically speaking, for any given $V_c < V_0 = \sum_{i=1}^n v_i$,  the   canonical penalty-duality method can produce global optimal solution to the bilevel minimization
problem $(\calP)$. However, if  $ \mu_c = V_c/V_0 \ll 1, $   any   iteration  method  could lead to unreasonable numerical solutions. In order to resolve this problem,
a volume reduction method  (VRM) method can be proposed:
\begin{verse}
(i) Introduce a volume reduction control parameter
$\mu \in (\mu_c,1)$ to produce a volume reduction sequence  $\{ V_{\gamma } = \mu V_{\gamma-1}\} $ ($\gamma =  1, \dots, \gamma_c$)
such that $V_{\gamma_c} = V_c$, $\calZ_\gamma = \{ \bz \in \{0,1\}^n| \;\;  \bv^T \brho \le V_\gamma\}$; \\
(ii) For a  given $V_\gamma \in [V_c, V_0]$ and  initial values $ (\bu^{\gamma-1}, \brho^{\gamma-1}) \in \calU_a \times \calZ_a$,
to find $(\bu^\gamma, \brho^\gamma)$ such that
 \end{verse}
 \begin{eqnarray*}
(\calP)^\gamma: \;\;    (\bu^\gamma, \brho^\gamma) &   = &    \arg \min \left\{ \PP (\bx,\brho) |\; \bx\in\calX_a, \;\; \bz \in \calZ_\gamma \right\}, \\
 s.t. \;\; \bx   & = & \arg \min \{\Pi(\bchi, \bz) | \;\; \bchi \in \calX_a \}.
\end{eqnarray*}

Generally  speaking,  for any given sequence $\{V_\gamma \}$ we should have
 \eb
 (\calP) = \lim_{\gamma\rightarrow \gamma_c} (\calP)^\gamma.
 \ee
Numerically, different volume sequence $\{V_\gamma\}$ may produce totally different numerical results as long as the decoupled  alternative iteration is used.
This is an  intrinsic difficulty for all   coupled bi-level  optimal design problems.
Based on the decoupled alternative iteration (DAI) and the volume reduction method (VRM),
 the canonical  duality  algorithm (CDT)  for solving the general bi-level optimization  problem $(\calP)$ can be proposed.

\begin{algorithm}
\caption{Canonical  Dual Algorithm for Bi-Level Knapsack Problem (CDT)}
\begin{algorithmic} [1]
\State  \textbf{Input parameters:}  $\mu$,     error allowances  $\omega  > 0$,
  $\brho^0 = \{1\} \in \real^n$,   and     $V_0 = \sum_{i=1}^n v_i  $. Let $\gamma = 1 $.
\State  \textbf{Solve}  the lower-level optimization problem  \eqref{eq-lowk}:
  \eb
  \bx^\gamma =  \arg \min \{ \Pi( \bu, \bz^{\gamma-1})   \;\; | \;\; \bx\in \calX_a \} .
  \ee
\State \textbf{Compute} $\bQ  = \bQ(\bx^{\gamma })$,  $\bb  = \bb(\bx^\gamma)$, and $V_{\gamma } = \max\{ V_c,  \mu V_{\gamma-1}\} $.

  \State \textbf{Solve}  the canonical  dual  problem $(\calP^d_q)^\gamma$ to obtain
   \eb
\bzeta^\gamma =    (\bsig^\gamma, \tau^\gamma) = \arg   \min \{  \PP^d_q (   \bzeta)    | \; \;\; \bzeta \in \calS^+_a \}
  \ee
     \State   \textbf{Compute} the upper-level  solution   $\bz^{\gamma}   $ by 
        \[
        \bz^{\gamma }   = \bG(\bsig^\gamma)^{-1} \btau(\bzeta^\gamma)
        \]
      \State    {\bf  If }   $| P_q (  \brho^{\gamma })  -  P_q(\brho^{\gamma-1})   | \le \omega $ and $V_\gamma  \le V_c$ ,
      then stop; Otherwise, continue.
                \State \textbf{Let}       $\gamma = \gamma+1 $. Go to Step 2.
\end{algorithmic} \label{alg-CDT}
\end{algorithm}
The canonical dual problem $(\calP^d_q)^\gamma$ in this algorithm could be either $(\calP^g_\alp)^\gamma$
or $(\calP^d_\beta)^\gamma$  for quadratic knapsack problem $(\calP_q)$.
For linear knapsack problem $(\calP_\ell)^\gamma$, its canonical dual could be either
 $(\calP^g_{\ell })^\gamma$ or $(\calP^d_{\ell\beta})^\gamma$.

\section{Application  to  Optimal Topology Design}\label{sec-example}
Topology design is the arrangement of the various elements (links, nodes, volumes, etc.) of a complex network or a discretized continuous system
  in multidisciplinary fields of
   communication, electronics,    optics,  structural, bio- and nano-mechanics (cf. \cite{bc,cml,gao-to17,kjy}).
 It was discovered recently by Gao   \cite{gao-to17} that the  optimal topology design  for  general  elastic structures  should be
  formulated as a bilevel knapsack problem:
   \begin{eqnarray}
   (\calP_{to}): &    & \min\{ \PP(\bx,\bz) = - \bb(\bx)^T \bz | \;\; \bx \in \calX_a, \;\;  \bz \in \calZ_a   \} \;\;\;\;\;\; \;\;\;\\
& &  s.t.  \;\; \bu \in \arg \min \left\{ \Pi( \bchi,  \bz)   \;\;   | \;\;  \bchi \in \calX_a  \right\} .\label{eq-lowto}
\end{eqnarray}
The leader design variable $\bz = \{ 0 , 1\}^n $ in this problem is the so-called density distribution such that if  $z_e = 1$, then  the $e$-th element of the elastic structure is  solid, otherwise, this element is void.
The follower  variable $\bx \in \calX_a \subset \real^m$ represents the  displacement vector, whose domain $\calX_a$ is a convex set, in which the   boundary condition is given (i.e. certain components of $\bx$ are fixed). The vector $\bb(\bx) =\{ \cc_e (\bx_e) \}   $ and its component  represents the strain energy   in the $e$-th element. 
For linear elastic structural, the total potential energy $\Pi(\bx,\bz)$ is a quadratic function of $\bx$:
\eb
\Pi(\bx,\bz) =  \half \bx^T \bK(\bz ) \bx - \bx^T \bff = \sum_{e=1}^n  \frac{z_e}{2} \bx^T_e \bK_e \bx_e   - \bx^T \bff,
\ee
 where $\bK_e$  and $\bx_e$ are  respectively  the stiffness matrix  and the nodal displacement vector of the $e$-th element;
 $\bv =\{ v_e \} \in \real^n_+ $ and $v_e \ge 0$ represents the volume of the $e$-th element;
$V_c$ is a desired volume.
Since $\bK_e $ is positive definite, we have $\bb(\bx) = \half  \{   \bx^T_e \bK_e \bx_e  \}^n \in \real^n_+ $ and the global stiffness matrix  $\bK (\bz) = \{ z_e \bK_e \} \in \real^{m\times m} $ is also positive definite for any given $\bz \in \calZ_a$.
The given vector $\bff \in \real^m$ is  the  external force.

Clearly, for a given displacement vector $\bx \in \calX_a$, the upper-level problem is a typical linear knapsack problem, which can be solved analytically
by the canonical duality theory.
While for a fixed design variable $\bz \in \calZ_a$, the lower-level minimization is the  well-known minimum total potential energy principle for linear elastic structures.
Since the total potential energy $\Pi(\bx,\bz)$ is a quadratic function of $\bx$, the lower-level solution can be given analytically by  $\bx = \bK(\bz)^{-1}  \bff$.
Thus, this  challenging topology design problem can be solved by the combination of the DAI-VRM and the canonical  duality method.
Since the linear knapsack problem can be solved analytically by Corollary 1 or 2, the computational complexity for solving
the lower-level solution  $\bx = \bK(\bz)^{-1}  \bff$ is about $O(m^3)$, the CDT is a polynomial time algorithm.

\begin{figure}[tb]
  \centering
\subfigure[2D beam]{\includegraphics[width=.45\textwidth, height=50mm]{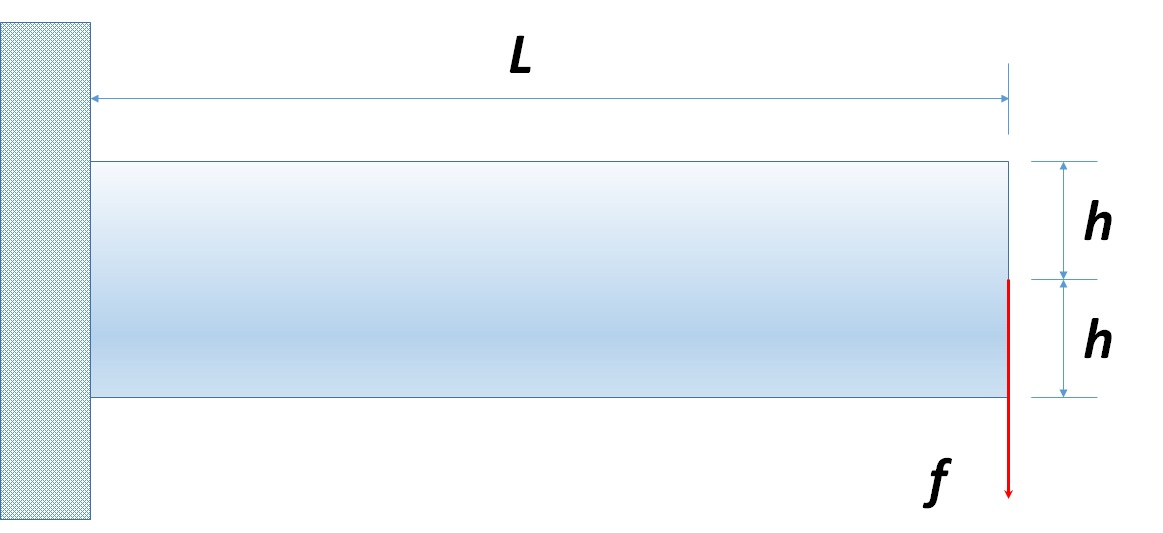}    }
\subfigure[3D beam]{\includegraphics[width=.45\textwidth, height=50mm] {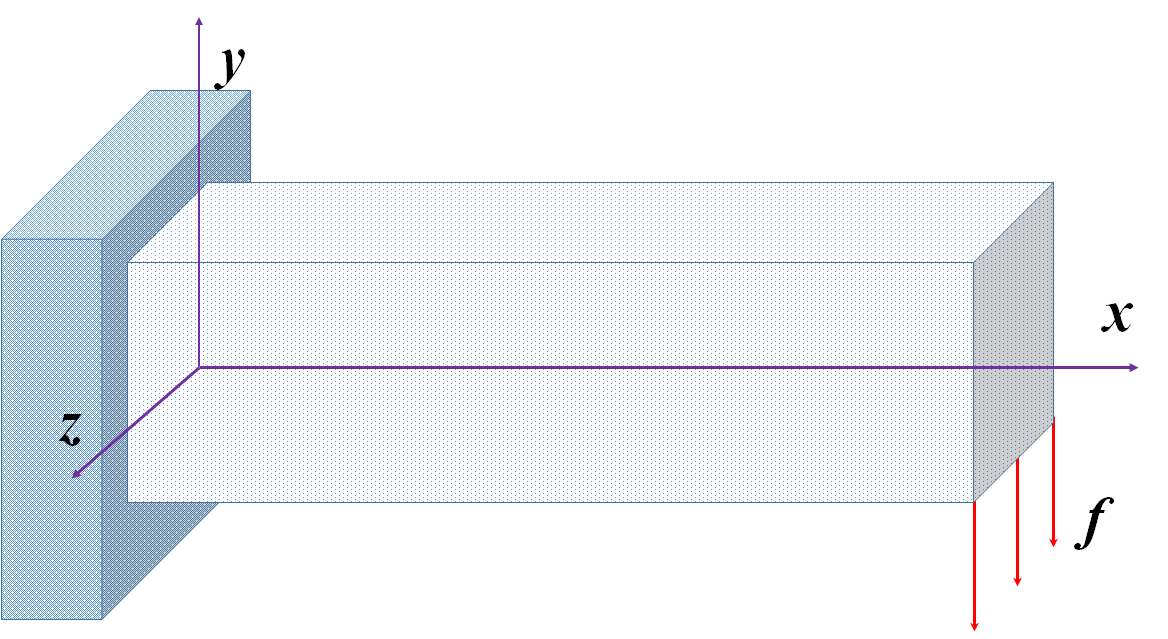}  }
     \caption{Design domains   for   long cantilever   beams with external load}  \label{fig-cant}
 \end{figure}

\begin{figure} [h]
  \centering
  \subfigure [mesh=$40 \times 10$,  $C=416.577,$ Time $=1.2165$sec] {
\includegraphics[width=.45\textwidth, height=30mm]{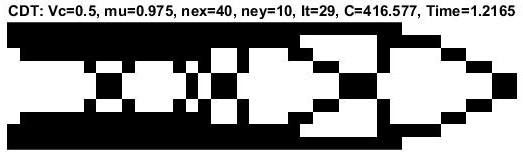}}
  \subfigure[mesh=$100 \times 30$,  $C=232.6$, Time$=2.96$sec] {
\includegraphics[width=.45\textwidth, height=30mm]{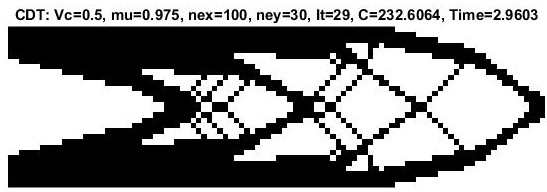}}
\subfigure[mesh=$180 \times 60$,   $C=171.77,$ Time $=7.69$sec ] {\includegraphics[width=.45\textwidth, height=30mm]{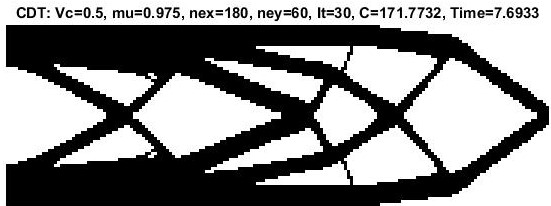}}
\caption{Optimal topology by different mesh resolutions with  $\mu_c = V_c/V_0 =  0.5$ and  $\mu=0.975$. Discrete variables $n=$mesh, continuous variables $m=2 (nex+1) (ney+1)$}
   \label{fig1}
\end{figure}

The proposed CDT algorithm  for topology design
 has been implemented in Matlab.
 The test examples are   the 2D and 3D   benchmark   cantilever beams    in optimal topology design of elastic structure
   (see  Fig.~\ref{fig-cant}).
Performance of the  CDT   method  is first  tasted   for different mesh resolutions.
Results in Fig. \ref{fig1} show    that for any given mesh resolutions,
 the CDT   method produces precise integer solutions.
Clearly,   the finer the resolution, the smaller  the stored energy $C=\bz^T \bc(\bx)$  with better material distributions.
This means that the fine structure has more capacity for external load.
For the mesh resolution  $180\times 60$ elements, we have $n=180\times 60 = 10,800$ discrete variables and $m=2\times (180+1) \times (60+1) = 21,960$ continuous variables, but the total computing time by a   HP labtop computer (with Processor Intel
Core I7-4810, CPU @ 2.80 GHz and memory 2.80 GB) is only $t = 7.69$ seconds.

\begin{figure} [h]
  \centering
   \subfigure[computational result by SIMP] {\includegraphics[width=.45\textwidth, height=40mm] {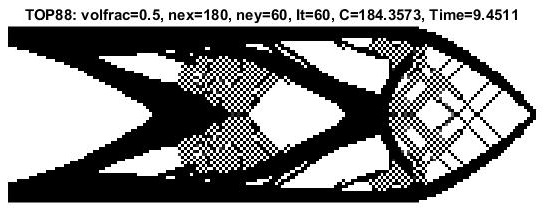} }
 \subfigure[computational result by BESO] {\includegraphics[width=.45\textwidth, height=40mm] {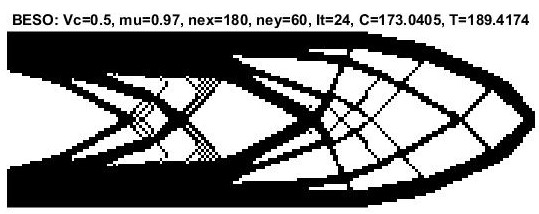} }
 \subfigure[Computational result by CDT]   {   \includegraphics[width=0.45\textwidth, height=40mm]{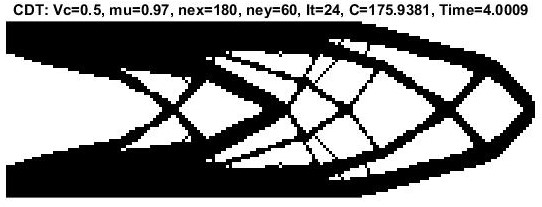} }
  \caption{2D beam structures by SIMP (a), BESO (b) and CDT (c) with  $n=180\times 60$, $m=2(180+1) \times (60+1)$.
The CDT is   twice faster than the SIMP, and about 47 times faster than BESO}
  \label{fig2}
\end{figure}

To compare with the commercial  codes  TOP88 (i.e. SIMP) and BESO \cite{Xie-Steven2}, we use the same mesh resolution of $180\times 60$.
The volume reduction rate $\mu=0.97$ for both BESO and CDT.
Computational results are reported in
Figure \ref{fig2}, which show clearly that the CDT method produces
geometrically simple and mechanically sound structure. 
Since the SIMP method is  a continuous relaxation approach, which cannot produce integer solutions,
the structure obtained by the SIMP code is ugly with a large area of grayscaling checkerboard  patterns.
By the fact that the BESO code is an evolutionary method for solving the linear knapsack problem by simply using comparison algorithm, although it can produce integer solution similar to that by CDT, it is not a polynomial time algorithm \cite{gao-to18},
therefore, it takes much longer time  (more than 30 times) than the CDT  method.

\begin{figure*}
\begin{tabularx}{0.5\textwidth}{ccc}
\hspace{-1.3em}
{
\begin{tabular}{lc}
\small  \subfigure[ CDT: $V_5=0.8587$, Time=4.97] {\includegraphics[scale=.35]{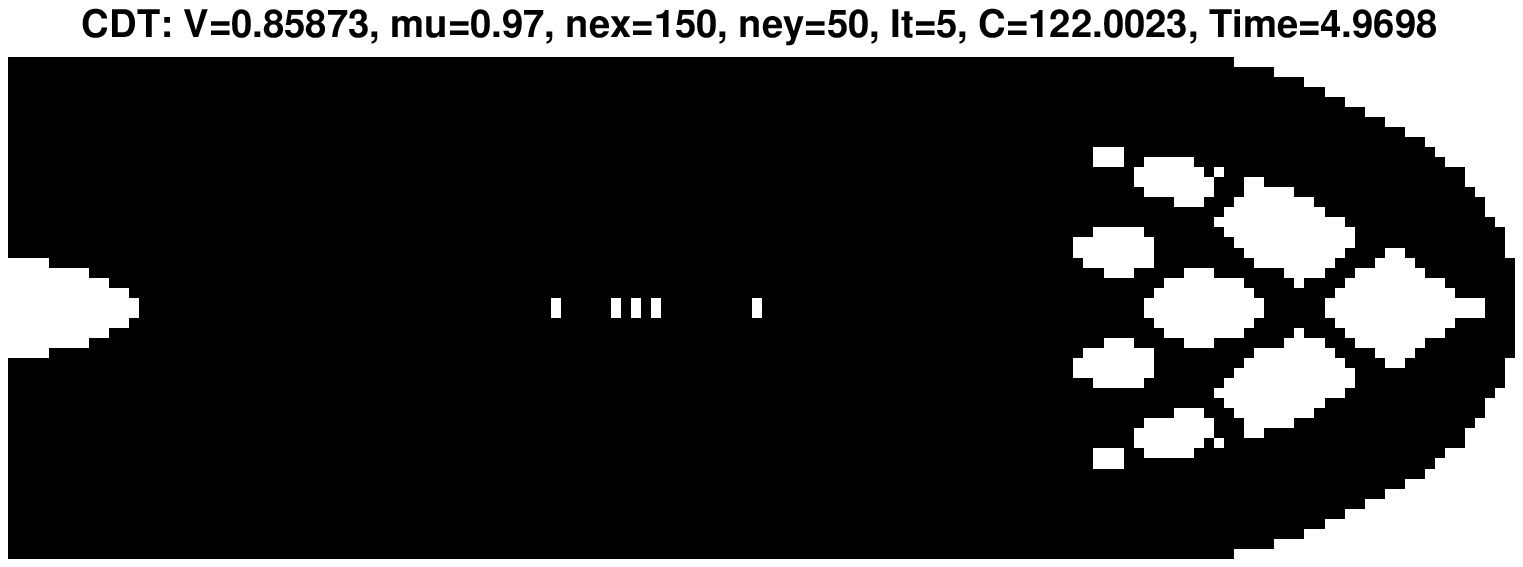}}\vspace{.4cm}\\
\small  \subfigure[ BESO: $V_5=0.8585$, Time=23.66] { \includegraphics[scale=.35]{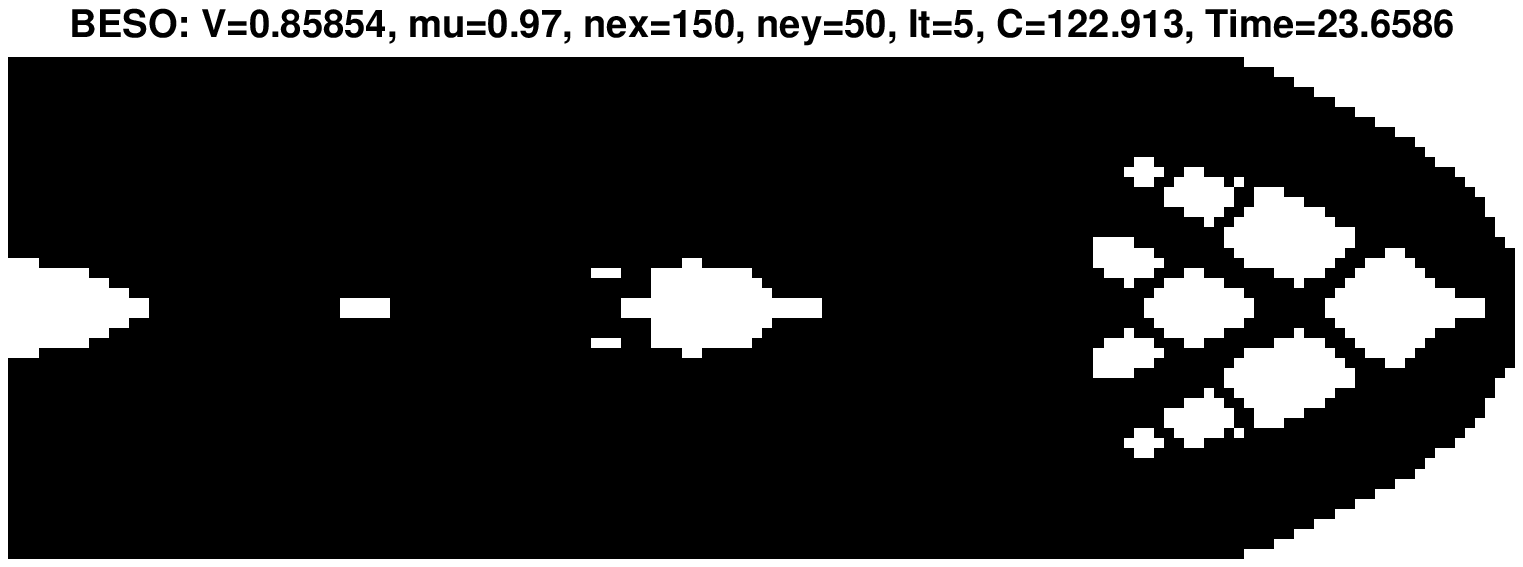}}\vspace{.3cm} \\
\small   \subfigure[SIMP:  Itration=5,  Time=3.07] {\includegraphics[scale=.35]{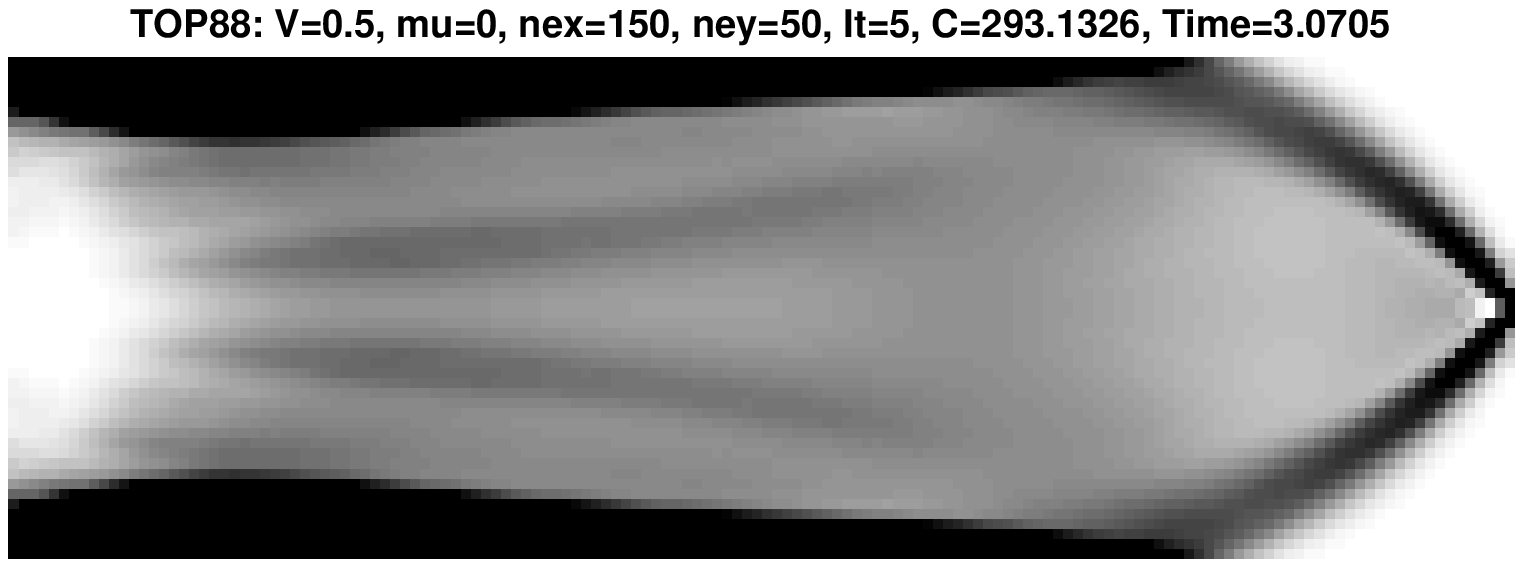}}
\end{tabular}}
&
\hspace{-2.5em}
{
\begin{tabular}{c}
\small \subfigure[ $V_{15}=0.63325$, Time=14.1421] {\includegraphics[scale=.35]{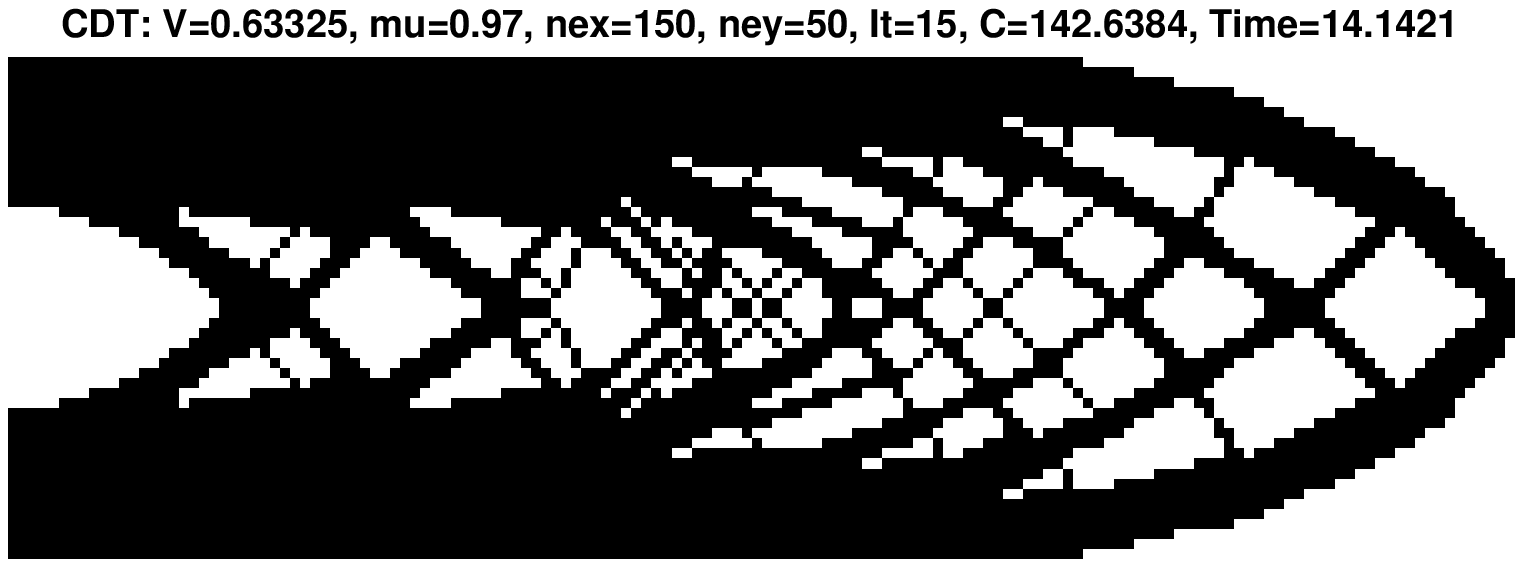}}\\
\small \subfigure[ $V_{15}=0.63343$, Time=73.812] {\includegraphics[scale=.35]{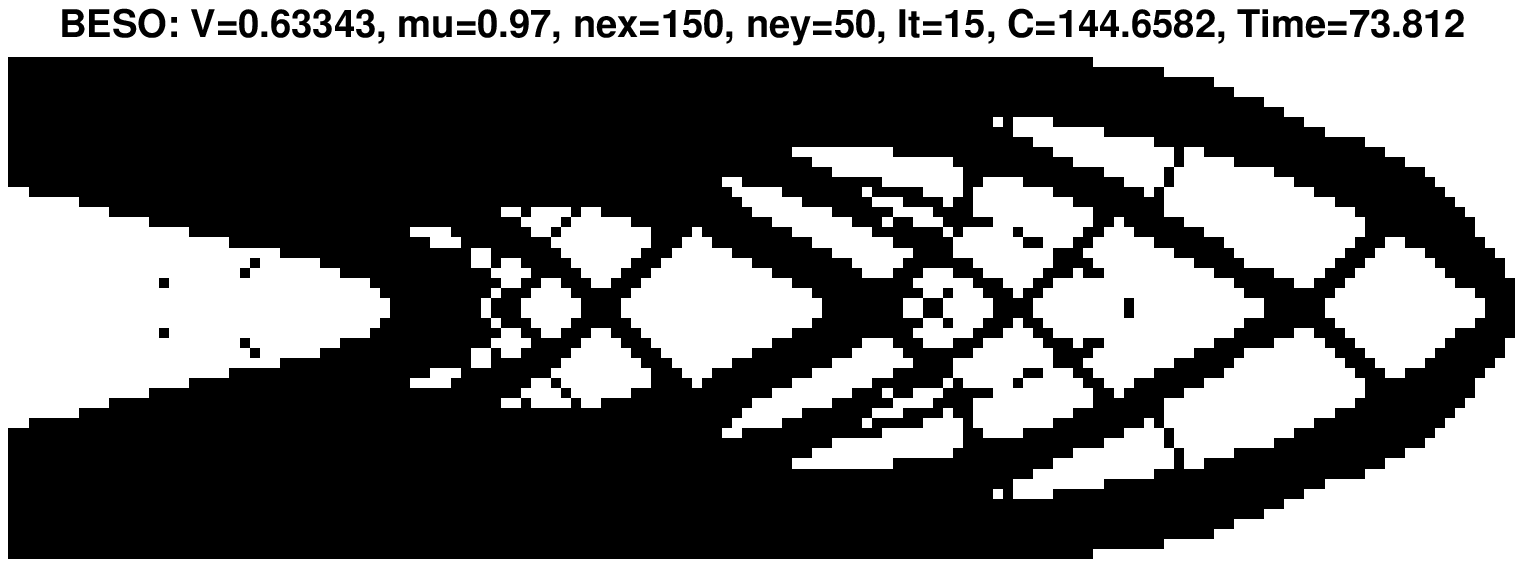}}\\
\small \subfigure[ Iteration=15, Time=11.1871] {\includegraphics[scale=.35]{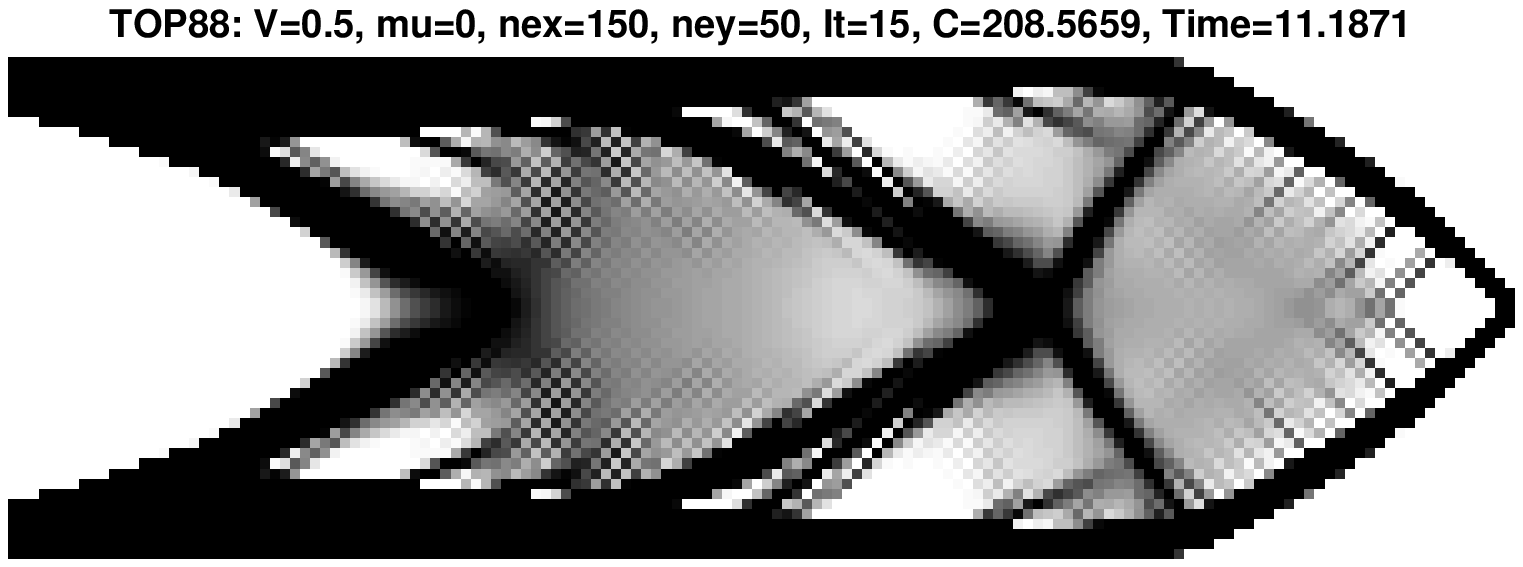}}
\end{tabular}}&
\hspace{-2.5em}
{
\begin{tabular}{c}
\small \subfigure[ $V_c=0.5$, Time=25.4657] {\includegraphics[scale=.35]{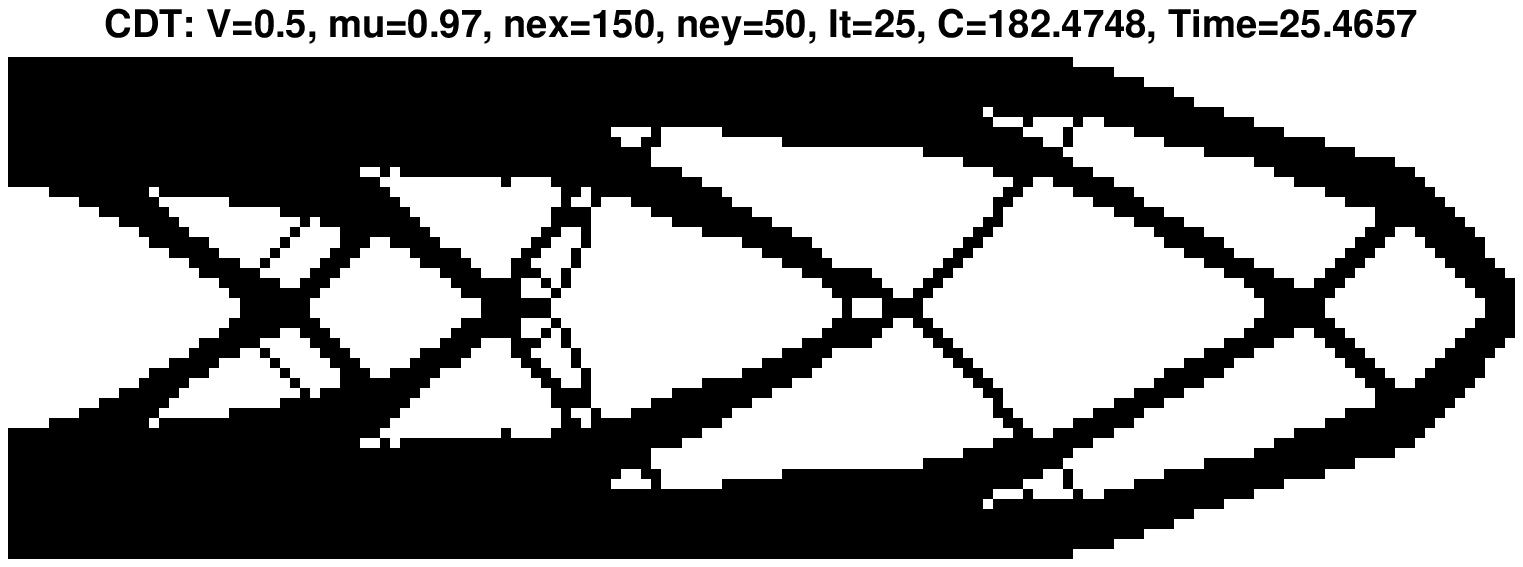}}\\
\small \subfigure[ $V_c=0.5$, Time=116.332] {\includegraphics[scale=.35]{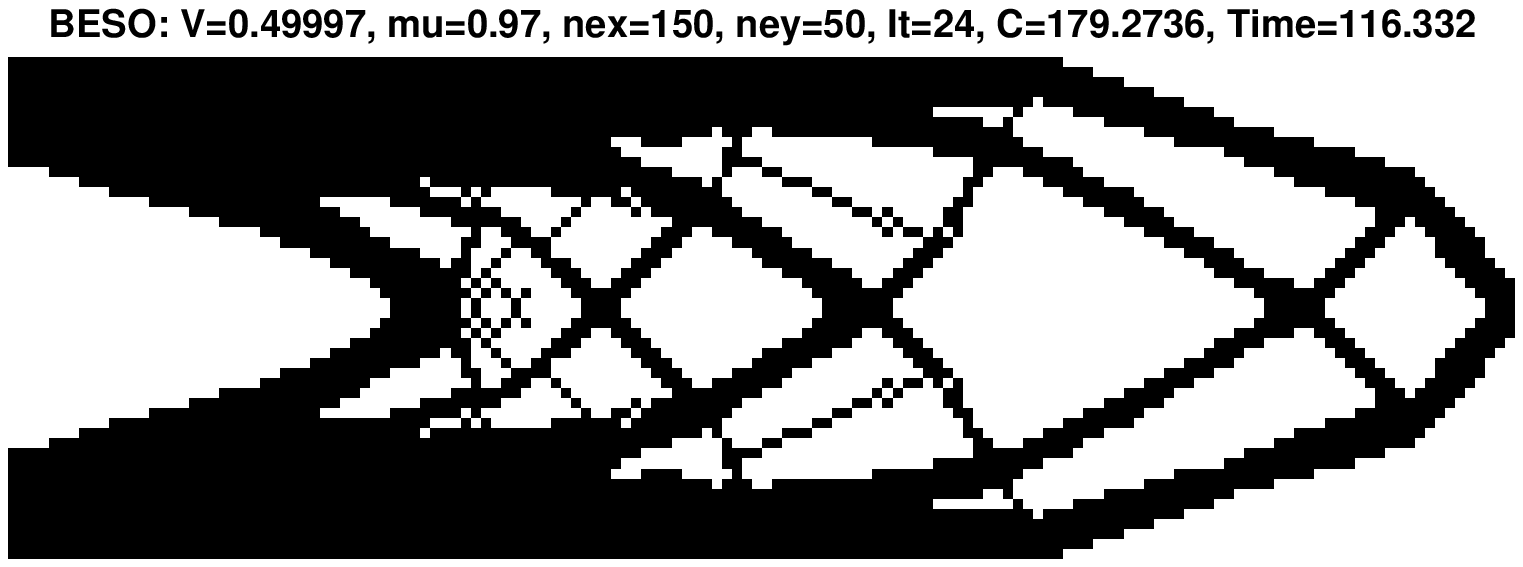}}\\
\small \subfigure[ Final Iteration, Time=104.0837] {\includegraphics[scale=.35]{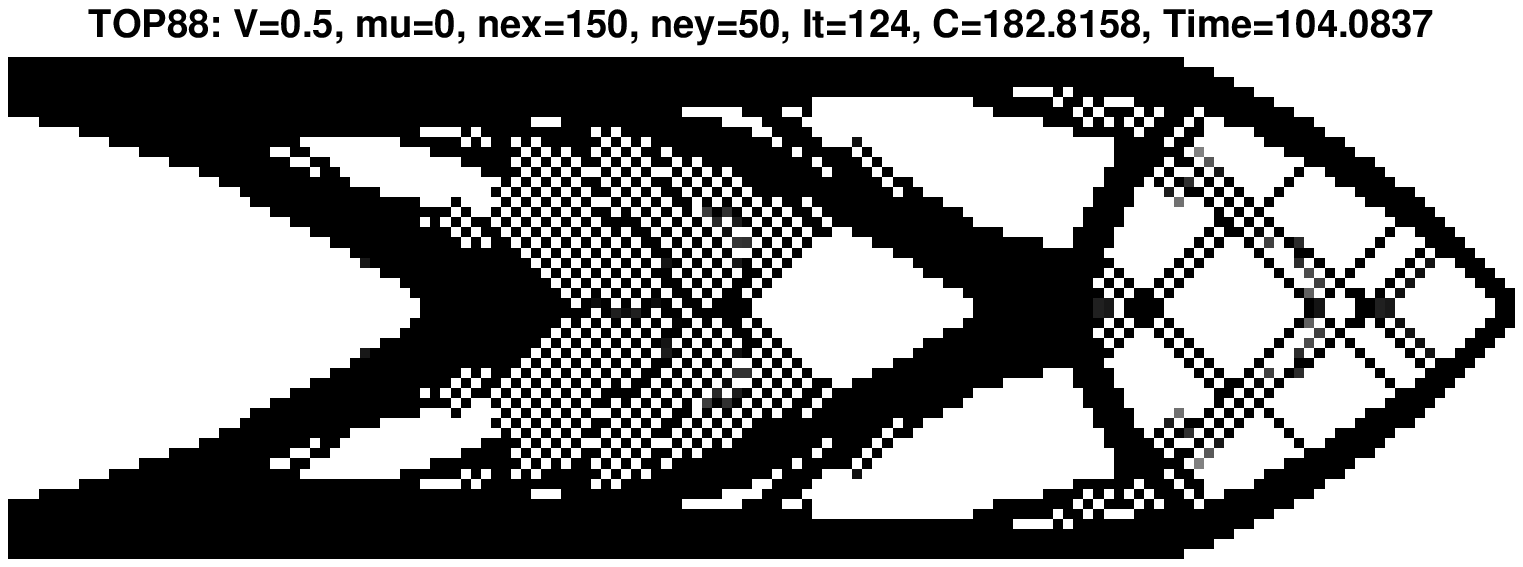}}
\end{tabular}}
\end{tabularx}
\caption{Topology evolutions for 2D  beam by CDT (top), BESO (middle), and SIMP (bottom)
with  $V_c=0.5$ and mesh resolution $150 \times 50$. The CDT is about four times faster than both BESO and SIMP.
The SIMP can produce only grayscale  topology since it is a continuous approach, i.e. $\bz \in (0, 1)^n$}\label{fig2d3}
\end{figure*}

In order to see the structures' evolution  during the iterations,
we tested the 2D beam with $V_c=0.5, \;\;  nex = 150, \;\; ney=50$.
Fig. \ref{fig2d3} shows the outputs of the structures by all methods  at the fifth, 15th and the final iterations.
Since the initial volume is $V_0 = 1$, the volume reduction rate is $\mu=0.97$ for both CDT and BESO,
therefore, the volumes at the 5-th and the 15-th iteration should be $V_5=0.859$ and $V_{15} = 0.633$, respectively.
 Results  in Fig. \ref{fig2d3} (a) and (b) show clearly that the structures produced by CDT algorithm are
final knapsack solutions if the desired knapsacks $V_c$  are $V_5 $ and $V_{15}$, respectively.
However, the structures produced by BESO are   broken with disconnected branches Fig. \ref{fig2d3}(b),
the structures by SIMP are uncertain with large areas of grayscales.

\begin{figure} [h]
  \centering
  \subfigure[SIMP: $C=178.688$] {
\includegraphics[width=.5\textwidth, height=40mm]{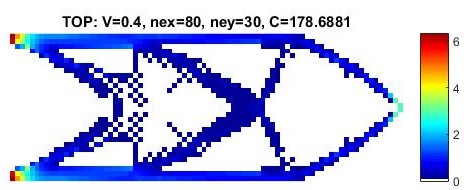}}
  \subfigure[BESO: $C=176.494$] {
\includegraphics[width=.5\textwidth, height=40mm]{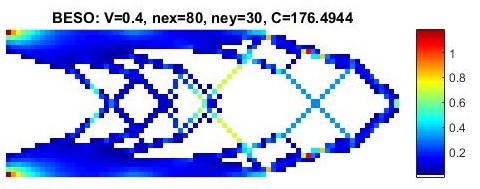}}
\subfigure[CDT: $C=170.869$] {\includegraphics[width=.5\textwidth, height=40mm]{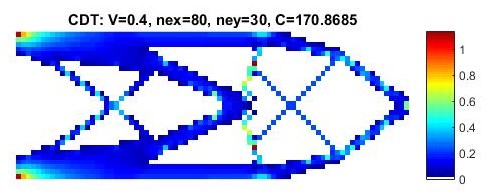}}
\caption{Structures and energy density distributions  by SIMP (a), BESO (b) and CDT (c) with $V_c=0.4$.}
  \label{fig23}
\end{figure}

 In order to look the   energy density  distribution $\bc = \{ c_e(\bu) \}$ in the optimal structures, we use the mesh resolution  $nex\times ney = 80 \times 30$.
 Fig. \ref{fig23} shows clearly that the CDT can produce mechanically sound structure with homogeneous distribution of strain energy density. Since the material distribution produced by SIMP is not mechanically reasonable, the overall strain energy stored in each element is about five times  higher than that by CDT (see the color bars in Fig. \ref{fig23}).  If we consider the dark-red level (= 1)  is the elastic limit, then the structure produced by SIMP is far beyond this limit.
This shows that the optimal structure by CDT has much potential to support even more  external load.
Although the BESO can produce better result than the SIMP, it can be used only for small-sized problems since it is not a polynomial time algorithm.

\begin{figure}
\subfigure[SIMP solution: C=4465.3]{\includegraphics[scale=0.45]{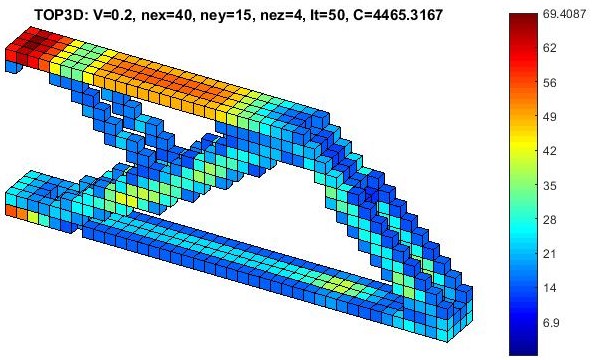}}\quad
\subfigure[SIMP solution: Front view]{\includegraphics[scale=0.45]{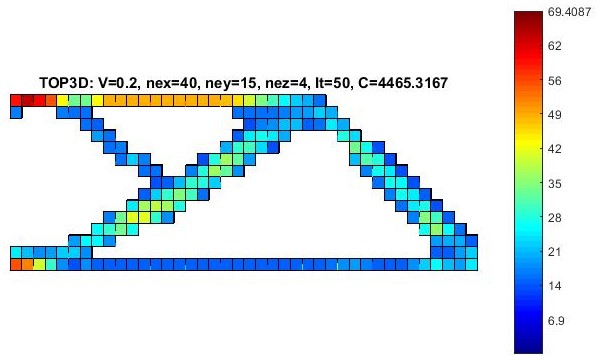}}
 \subfigure[CDT solution: C=2792]{\includegraphics[scale=0.45]{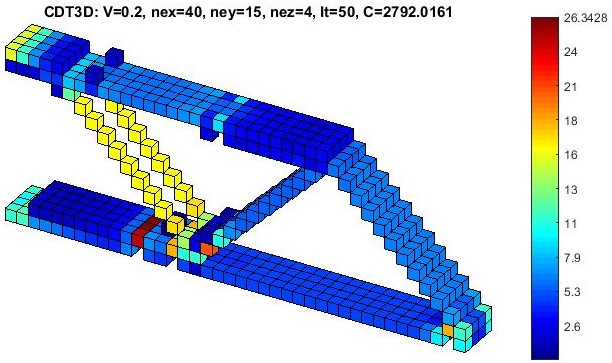}}\quad
\subfigure[Front view of CDT solution]{\includegraphics[scale=0.4]{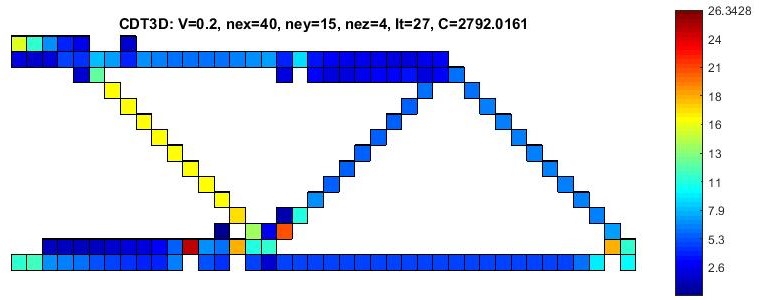}}
\caption{Structures  for 3-D   beam with mesh=$40\times 15 \times 4$   and $V_c =0.2$. The SIMP beam is overstaffed with  many red elements at scale=69.4 (highly stressed).
The CDT beam is geometrically elegant with very few red elements at scale=26.34, so it can support more external load. }\label{fig3D1}
\end{figure}

For the 3D beam,  we only compare the CDT with the SIMP since the BESO is computationally expensive.
First, we use the mesh resolution  $40\times 15 \times 4$ and volume fraction  $V_c = 0.2$.
Fig. \ref{fig3D1} shows that the optimal topology produced by CDT is not only geometrically elegant, but also mechanically sound with homogeneous distribution of energy density and very few red elements. From the color bar we can see that the
 red-level is only at scale of 26.34. However, the topology produced by SIMP is  overstaffed with many red elements at scale of 69.4.   This shows that the structure produced by the CDT can support external load three times  than the structure by SIMP.

\begin{figure}[tb]
\centering
\subfigure[CDT solution: C=13192.7, Time=1.56min]{
\includegraphics[width=.45\textwidth] {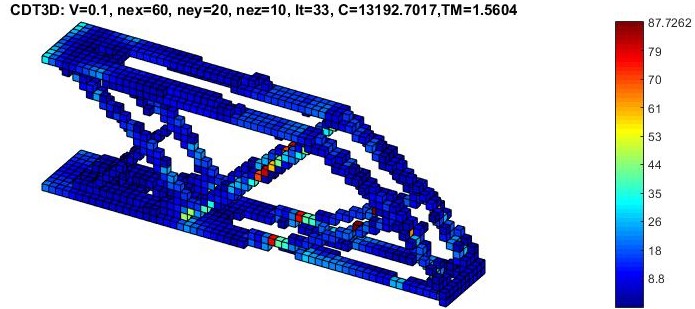} }\quad
\subfigure[Front view of CDT solution ]{
\includegraphics[width=0.4\textwidth, height=35mm] {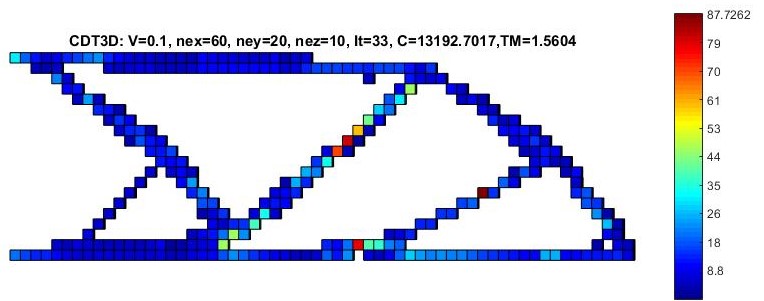} }
\qquad
\subfigure[SIMP solution at the $50$-th iteration: C=21960, Time=12.2min]{
\includegraphics[width=.45\textwidth] {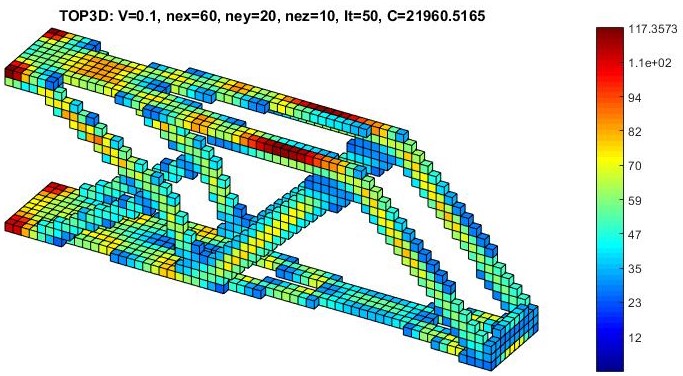} }\quad
\subfigure[Front view of SIMP solution ]{
\includegraphics[width=0.4\textwidth, height=45mm] {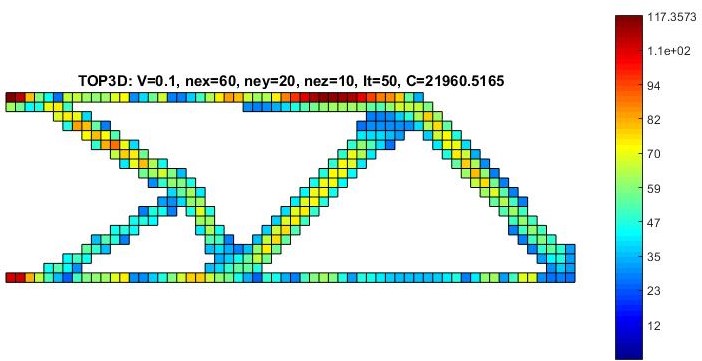} }
\caption{3D beam structures  by CDT and SIMP with $V_c=0.1$
and mesh resolution  $ 60 \times 20 \times 10$. Since the strain energy stored in the CDT beam ($C=13192.7$)  is only
about a  half of that stored in the SIMP beam ($C=21960$), the CDT beam
 has much more potential to support  more external load.}
\label{fig4}
\end{figure}

\begin{figure}[tb]
\centering
\subfigure[SIMP solution ($It=199$): C=21426, Time=19.8min]{
\includegraphics[width=.40\textwidth] {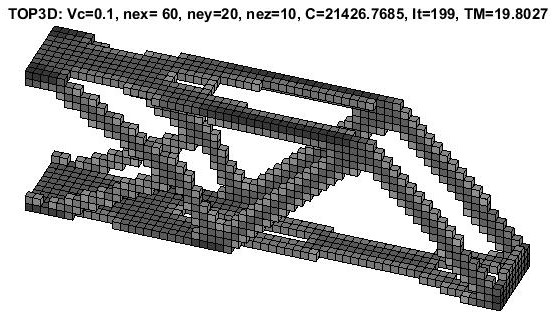} }\quad
\subfigure[Front view of SIMP solution ]{
\includegraphics[width=0.4\textwidth, height=35mm] {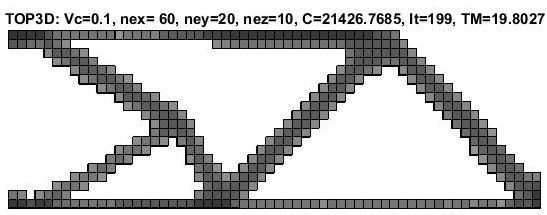} }
\caption{3D beam structures  by  SIMP with $V_c=0.1$
and mesh resolution  $ 60 \times 20 \times 10$.  }
\label{fig5}
\end{figure}

Now we increase the mesh resolution to
   $60\times 20 \times 10$ and reduce the volume fraction to $V_c = 0.1$.
  Fig. \ref{fig4} shows results   obtained by  the commercial code TOP3D (SIMP)  and CDT3D     (with $\mu=0.93$).
Since  the  0-1 density distribution produced  by the  CDT is much more mechanically reasonable than the one by SIMP,
most of elements in Fig. \ref{fig4}(a)  are in blue color with only a few red elements at the scale= 88.
The beam by   SIMP has more red  elements  at the scale =117.
The structure shown in Fig. \ref{fig4}(c) is the out put of TOP3D at the  $It=50$-th iteration, which is almost the same as the final result
($It=199$) shown in Fig. \ref{fig4}(e). This result shows that the SIMP has a very slow rate of  convergence.
Also the CDT is more than 10 times faster than the TOP3D.
The uncolored Fig. \ref{fig5} shows  that most   elements  are grayscale (i.e. $0< z_i < 1$)  since the SIMP is a continuous approximation.

Detailed study on the canonical duality theory for solving topology design  problems and a 66-line  Matlab code  of the CDT algorithm   will be  given   in
\cite{gao-r-l-to}.

\section{Concluding remarks  and open problems}\label{sec-con}
Based on a  combination of an alternative iteration, a volume reduction and    the novel canonical duality theory,
a new powerful   method  is proposed   for solving  bilevel optimization  governed by knapsack problems.
The main theoretical contributions are presented in Sections  3 and 4, i.e.  by using  the canonical  duality method,  the well-known knapsack problems can be  solved
 analytically in terms of their canonical dual solutions.
The existence and uniqueness of this general analytical solution are proved for both quadratic and linear knapsack problems.
It shows that for any given $\bb,\;\bv\in \real^n$, if there exists a canonical dual solution  $\tau_c >0$ such that $\tau_c \bb \neq \bv$, the  linear knapsack problem
can be solved  uniquely in polynomial time.
Therefore, a polynomial  algorithm CDT  is proposed for solving the bilevel knapsack problems.
Its novelty is demonstrated by solving a benchmark  problem  in topology design.
Numerical results indicated   that the volume reduction method is essential for solving bilevel optimization  governed by knapsack problems.
Both theoretical and numerical results verified  that the   CDT  is a powerful method for studying
general bilevel  optimization problems.

Due to the intrinsic  coupling effect in  the  bilevel  optimization  problem $(\calP)$,
  although  for each  given lower level variable $\bx_k\in \calX_a$,  the  upper level  knapsack problem  $\bz_k = \arg \min \{ P(\bx_k, \bz) | \; \; \bz   \in \calZ_a  \} $
  can be solved analytically by the canonical duality theory, the solution sequence $\{\bx_k,\bz_k\}$  may not converge to  the global optimal solution of  $(\calP)$
  since the decoupled alternative iteration method is adopted
 in Algorithm 1.
This is the main reason that  the numerical solutions produced by the CDT  depend  sensitively on  the  volume reduction  rate $\mu$.
 This is an open problem and deserves   theoretical study in the future.

\section*{Acknowledgement}
The author would like to sincerely acknowledge the important comments and suggestions from  the editor and anonymous
reviewers.
This paper is based on several  invited/keynote/plenary
lectures presented at international conferences including
the Symposia on Intelligent Technologies for Advancing and Safeguarding Australia, 15 August 2017, Deakin University,
and the  Second In. Conf.  on Modern Mathematical Methods and High Performance Computing in Science and Technology, 4-6 January, 2018, New Delhi, India.
 Organizers' hospitality and financial supports are sincerely acknowledged.
 This research is supported by US Air Force Office for Scientific Research (AFOSR)  under the grants   FA2386-16-1-4082 and FA9550-17-1-0151.

\begin{IEEEbiography}[{\includegraphics[scale=.18]{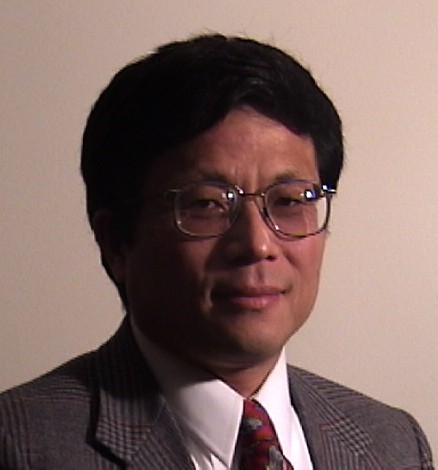}}]{David Yang Gao}
is the Alex Rubinov Professor of Mathematics at the Federation University  Australia. After he received his PhD from Tsinghua University in Engineering Mechanics and Applied Mathematics, he has held research and teaching positions in different institutes including MIT, Yale, Harvard, the University of Michigan, and Virginia Tech.
His research has been mainly focused on modeling and fundamental principles of duality and triality in multidisciplinary 
fields of engineering science, mathematics, physics, information technology, natural and social 
complex systems. He has published about 15 books and  200 scientific and philosophic papers 
(about  50\% are single authored). 
 His main research contributions include a canonical duality-triality theory, which can be used to model complex phenomena within a unified framework. This theory has been used successfully for solving a large class of nonconvex/nonsmooth problems in nonlinear
sciences as well as a series of well-known NP-hard problems in global optimization and computational science. The main part of the canonical duality theory, i.e., the complementary-dual variational principle he proposed in 1997, solved a 50-year open problem in nonlinear elasticity and  now is called the Gao principle in literature. This principle plays an important role in nonconvex analysis and computational mechanics. 
Professor Gao is a founding editor for Springer book series of Advances in Mechanics and Mathematics. He has served as an associate editor for several international journals of applied math, optimization, solid mechanics, dynamical systems, industrial and management engineering. He  is a founding   vice president and secretary general  of the International Society of Global Optimization.
\end{IEEEbiography}
\end{document}